∫%
\documentclass[journal]{IEEEtran}

\usepackage{algorithm}
\usepackage{algorithmic}
\usepackage{xcolor,soul,framed} 

\colorlet{shadecolor}{yellow}
\usepackage{color}
\usepackage[pdftex]{graphicx}
\graphicspath{{../pdf/}{../jpeg/}}
\DeclareGraphicsExtensions{.pdf,.jpeg,.png}
\usepackage{mathrsfs}
\usepackage[cmex10]{amsmath}
\usepackage{amsmath}
\usepackage{array}
\usepackage{mdwmath}
\usepackage{mdwtab}
\usepackage{eqparbox}
\usepackage{url}
\usepackage[yyyymmdd,hhmmss]{datetime}
\usepackage{multirow}
\allowdisplaybreaks

\usepackage{amssymb} 
\usepackage{verbatim}
\usepackage[noadjust]{cite}

\usepackage{makecell} 

\usepackage{pifont}
\newcommand{\xmark}{\ding{53}}%

\usepackage{nomencl}
\makenomenclature

\newcommand{\y}[0]{\boldsymbol{\textcolor{blue}{y}}}
\newcommand{\s}[0]{\boldsymbol{\textcolor{blue}{s}}}
\newcommand{\DeltaP}[0]{\boldsymbol{\textcolor{blue}{\Delta P}}}

\newcommand{\x}[0]{\boldsymbol{\textcolor{red}{x}}}
\newcommand{\rtau}[0]{\boldsymbol{\textcolor{red}{\tau}}}
\newcommand{\rb}[0]{\boldsymbol{\textcolor{red}{b}}}
\newcommand{\q}[0]{\boldsymbol{\textcolor{red}{q}}}

\usepackage{ifthen}
  \renewcommand{\nomgroup}[1]{%
  \item[\bfseries
  \ifthenelse{\equal{#1}{D}}{Definitions/Abbreviations}{%
  \ifthenelse{\equal{#1}{S}}{Indices/Sets}{%
  \ifthenelse{\equal{#1}{P}}{Parameters}{%
  \ifthenelse{\equal{#1}{V}}{Variables}{}}}}%
  ]}

\usepackage{threeparttable}  

\usepackage{flushend} 

\makeatletter
\newcommand\fs@norules{\def\@fs@cfont{\bfseries}\let\@fs@capt\floatc@ruled
  \def\@fs@pre{}%
  \def\@fs@post{}%
  \def\@fs@mid{\kern3pt}%
  \let\@fs@iftopcapt\iftrue}
\makeatother
\floatstyle{ruled}
\restylefloat{algorithm}

\usepackage{titlesec}
\usepackage{hyperref}

\makeatletter
\renewcommand\paragraph{\@startsection{paragraph}{5}{\z@}%
  {3.25ex \@plus1ex \@minus.2ex}%
  {-1em}%
  {\normalfont\normalsize\bfseries}}
\def\toclevel@paragraph{5}
\def\toclevel@paragraph{6}
\def\l@paragraph{\@dottedtocline{5}{10em}{5em}}
\makeatother

\usepackage{fancyhdr}

\fancypagestyle{firstpage}{%
  \fancyhf{} 
  \fancyfoot[C]{\tiny \color{blue} This work has been submitted to the IEEE for possible publication. Copyright may be transferred without notice, after which this version may no longer be accessible.\\ © 2021  IEEE.  Personal use of this material is permitted.  Permission from IEEE must be obtained for all other uses, in any current or future media, including reprinting/republishing this material for advertising or promotional purposes, creating new collective works, for resale or redistribution to servers or lists, or reuse of any copyrighted component of this work in other works}
  
}
\begin{document}
\bstctlcite{IEEEexample:BSTcontrol}
    \title{Joint Design for Electric Fleet Operator and Charging Service Provider:  Understanding the Non-Cooperative Nature}
    %
    
    
    %
    
  \author{
  Yiqi~Zhao$^\ast$, 
  Teng Zeng$^\ast$, 
  Zaid Allybokus,
  Ye~Guo, \IEEEmembership{Senior Member,~IEEE,}
  Scott Moura, \IEEEmembership{Member,~IEEE}
      
   \thanks{$^\ast$ equal contributions. This work was supported in part by Total S.E. \textit{Corresponding author: Teng Zeng (tengzeng@berkeley.edu)}.}
   \thanks{T. Zeng and S. Moura are with the Department of Civil and Environmental Engineering, University of California, Berkeley 94720 CA.}
   \thanks{Y. Zhao, Y. Guo and S. Moura are with the Smart Grid and Renewable Energy Laboratory, Tsinghua-Berkeley Shenzhen Institute, Shenzhen, China.}
   \thanks{Zaid Allybokus is with Total S.E. Smart Mobility group.}
}


\maketitle
\thispagestyle{firstpage}
\begin{abstract}

This work proposes a new modeling framework for jointly optimizing the charging network design and the logistic mobility planning for an electric vehicle fleet. Existing literature commonly assumes the existence of a single entity -- the social planner, as a powerful decision maker who manages all resources. However, this is often not the case in practice. Instead of making this assumption, we specifically examine the innate non-cooperative nature of two different entities involved in the planning problem. Namely, they are the charging service provider (CSP) and the fleet operator (FO). To address the strategic interaction between entities, a bi-level mixed integer program is formulated, with the CSP/FO's problem expressed in the upper/lower levels respectively, in a joint decision making process. These decisions involve the CSP's infrastructure siting, sizing, substation capacity upgrades, the FO's fleet composition, vehicle routing, charging, and delivery assignment. To relieve computational burdens, we utilize a double-loop solution architecture to iteratively reach optimality. We conduct detailed numerical studies on a synthesized small network and the simulation results reveal the unique aspects of this two-entity framework.
This modeling perspective can be generalized to other system design problems with two interacting agents planning and operating resources across networks.


\end{abstract}


%
\IEEEpeerreviewmaketitle

\section{Introduction}
\IEEEPARstart{D}{ecarbonization} of the transportation sector is an important step towards alleviating climate change. In the U.S., about 28\% of the total greenhouse gas emissions is contributed by transportation \cite{useap_trans}. According to the California Air Resources Board (CARB), commercial trucks are responsible for 80\% of the diesel soot emitted, leading the largest source of air pollution from vehicles \cite{carb_truck}. Hence, a significant step to cut emissions is electric commercial vehicles, specifically E-trucks, as part of a sustainable supply chain. As a result, CARB has voted to rule out the sales of any fossil fuel trucks by 2045 and to force truck makers to begin the transition in 2024 \cite{carb_truck}. 

Along with enforcing policy orders, many logistic and delivery companies (we refer to as fleet operators, FOs) and Charging Service Providers (CSPs) are committed to transportation electrification. To realize profit maximization, it is more important to have effective communications between these two entities. The CSPs, with knowledge of spatial-temporal charging demands, could strategically construct their charging network to accommodate the needs; whereas the FOs, whose electrified trucks are constrained by driving range,  would consider charging en-route but with minimal detours. 

\subsection{Literature Review}\label{subsec:lit_review}
In this section, we are going to review a body of literature that we have identified as not only relevant but also crucial to understand our problem. On one hand, to consider commercial E-trucks routing, a portfolio of attributes can be considered, including \textit{homogeneous/heterogeneous fleet, range, partial/fully charging time, delivery time windows, etc}. On the other hand, to consider charging infrastructure planning, another set of attributes are considered, such as \textit{the station locations, power constraints, etc}. The two entities, FO and CSP, are entangled through charging events and an extensive body of literature (\cite{yang2015battery, hof2017solving, schiffer2016ecvs, schiffer2018designing, schiffer2017adaptive, schiffer2018strategic, schiffer2017electric, li2015multiple, hiermann2016electric}) has accounted for this interactions. In the field of Operation Research and Electrical Engineering, this is called the electric location routing problem (E-LRP), an extension to pure electric vehicle routing problem (E-VRP). For the classical vehicle/location routing problems, we refer interested readers to these two survey papers \cite{prodhon2014survey, schiffer2019vehicle}.

In the aforementioned literature, each work varies focus slightly and considers a subset of the above entity-specific attributes. J. Yang and H. Sun \cite{yang2015battery} were the first to investigate the E-LRP, where the location of battery swapping stations (BSS) was jointly optimized together with the routing of a homogeneous E-trucks fleet. The computational results of the work were later improved by Hof, Schneider and Goeke \cite{hof2017solving}. M. Schiffer et al. conducted a series of research on E-LRPs. Each publication in this series has a different focus. For example, \cite{schiffer2016ecvs} incorporated real-world data to address the competitiveness between E-trucks and ICEVs,  \cite{schiffer2018designing, schiffer2017adaptive} considered deployment of multiple types of facilities (replenishment, recharging, and combined type facilities), \cite{schiffer2018strategic} addressed uncertainty using robust optimization , \cite{schiffer2017electric} used different planning objectives. Authors in \cite{li2015multiple} further considered multiple types of charging facilities in E-LRP with time windows (E-LRPTW). Paper \cite{hiermann2016electric}, on the other hand, investigated the effects of heterogeneous fleets on a similar E-VRPTW setting with a full recharge scheme. These works inevitably assume the existence of a powerful social planner, who is capable of coordinating all the tasks. However, this is often not possible in practice. Instead, the FO and the CSP are more likely to be separate organizations with misaligned incentives, leading to non-cooperative behavior. In this work, we specifically capture these dynamics and model it as a leader-follower Stackelberg Game, which to the authors' best knowledge, has never been studied. 

Furthermore, the modeling approaches in the above works closely resemble each other and are the natural stems from the classic traveling salesman problem (TSP). The abstracted network is often called the customer-node based network, where customer nodes are the graph representatives and constrained to be visited once and exactly once. Additional features like range limits and charging speed for E-trucks are easily incorporated via supplemental constraints. 

On the other hand, while the customer-node based network is classic and easy to adopt, the shortfall is prominent - lack of flexibility in tracking temporal events, such as charging. Since every node is associated with one specific set of entry and exit times for one vehicle, the temporal sense of simultaneous visits or queuing at a charging station node is dismissed. Adding trackers, e.g. indicator functions, is inevitable to address this issue. However, this makes the problem highly nonlinear and hence the solution quality cannot be guaranteed. Alternatively, references \cite{mahmoudi2016finding, lu2016solving} adopted the idea of layered graphs and proposed state-space-time/resource-space-time expanded networks to embed discretized resource values (energy consumption, time, etc.) when defining nodes. In this case, resource constraints are directly encoded in the expanded network model and time-dependent consumption patterns can be characterized. However, such modeling flexibility comes at a cost of significantly increased network size, and subsequently the problem scale. To plan a charging network, authors from \cite{zhang2019joint} took a macroscopic point of view with traffic flow and designed another way to expand the transportation network. In this network, all reachable nodes are extended with hyper-arcs to model EVs' feasible routes. Then, a joint fleet sizing and charging system planning method for autonomous electric vehicles was proposed. However, the core, which is the extended network, requires full recharge and predefined homogeneous battery capacities. Recently, the authors of \cite{kancharla2020electric} introduced a novel mixed integer linear programming model for the E-VRP with load-dependent charging patterns. The proposed formulation allows multiple visits of charging stations without expanding the network into higher dimensions, thus helping to reduce the problem scale. Although this approach neatly relaxes the aforementioned restriction, it cannot be directly applied to a setting where locations of charging facilities are unknown. To conclude, though charging station location planning for E-trucks has been studied, incorporating station size and capacity upgrade remains as gap in this field of research.

We have reviewed a series of literature and identified the remaining gaps in the community. In Table \ref{tbl:aspects_summary}, we summarize the aspects covered by some representative works and compare with ours. 



\begin{table}[!t]
\renewcommand{\arraystretch}{1.3}
\caption{Overview of considered aspects in existing literature}
\label{tbl:aspects_summary}
\centering
\begin{threeparttable}
\begin{tabular}{|c|ccccccc|}
\hline
Ref. & \makecell{View-\\point} & \makecell{Fleet\\Design} & \makecell{CS\tnote{*}\\Design} & \makecell{Charging\\Option} & \makecell{Logistic \\Cons.} & \makecell{TW\tnote{*}} &\makecell{CS\\Cons.}\\
\hline
\cite{hiermann2016electric} & FO\tnote{*} & heter & \xmark  & Fully & \checkmark & \checkmark & \xmark\\
\hline
\cite{kancharla2020electric} & FO & \xmark & \xmark & Partial & \xmark  & \checkmark & \xmark\\
\hline
\cite{bruglieri2019green} & FO & \xmark & \xmark  & Fully & \xmark & \checkmark & \checkmark\\
\hline
\cite{schiffer2018strategic} & SP\tnote{*} & homo & site & Partial & \checkmark & \checkmark & \xmark\\
\hline
\cite{yang2015battery} & SP & \xmark & site & Fully & \checkmark & \xmark& \xmark\\
\hline
\cite{lu2016solving} & SP & \xmark & site & Partial & \checkmark & \checkmark & \xmark \\
\hline
\cite{zhang2019joint} & SP & homo & size & Fully & \xmark & \xmark & \checkmark\\
\hline
Our & \makecell{Multi-\\players} & heter & \makecell{site\\+size} & Partial & \checkmark & \checkmark & \checkmark\\
\hline
\end{tabular}
\begin{tablenotes}
\footnotesize
\item[*] CS: charging station; TW: time window at customer points; FO: fleet operator; SP: social planner.
\end{tablenotes}
\end{threeparttable}
\end{table}

\subsection{Contributions}
This paper proposes a novel framework to jointly consider the charging infrastructure deployment and E-truck fleet design, as well as routing strategies. While the comprehensive settings will be simultaneously decided by the framework, the sense of non-cooperation and contest is maintained to clearly distinguish the two players: the CSP and the FO. Comparing to existing literature, the major innovations of this paper are:

\begin{enumerate}
  \item Instead of assuming the existence of a powerful single entity who owns both the fleet and the charging network, a two-party model with the charging service provider and fleet operator is adopted. They have their own objectives and their interactions are captured via a Stackelberg game, whose results are closely analyzed. The necessity of such modeling perspective is revealed.

  \item We propose an innovative partial time expanded network (PTEN) model on top of the customer-node based network \cite{yang2015battery, hof2017solving, schiffer2016ecvs, schiffer2018designing,schiffer2017adaptive, schiffer2018strategic, schiffer2017electric, erdougan2012green, li2015multiple, hiermann2016electric}. This enables us to track the simultaneous charging activities of E-trucks at each location without introducing any nonlinearity. The overall model is thus kept in the domain of mixed integer linear programming, which attracts computation advantages.
  \item The proposed model not only decides locations of the charging stations but also their sizes. Based on that, the upgrade cost of transformers is also incorporated in the CSP cost calculation, which is actually an important factor in real-world operation but has been neglected in past research.

  \item A double-loop solution framework is designed accordingly. The outer loop adopts the idea from \cite{zeng2014solving} to capture the dynamics of the Stackelberg game while ensuring convergence. The inner loop reformulates the fleet design and location/routing problem as generic set-partitioning problems, which are solved using column generation together with customized labeling algorithm. 

\end{enumerate}

The remainder of this paper is organized as follows: Section \ref{sec:prob_def_sys_model} gives formal problem definition and the proposed partial time expanded network modeling method. Based on the system model, Section \ref{sec:math_model_bi_level} presents the detailed mathematical formulations of the planning problem. The solution algorithm for the formulated model is then described in Section \ref{sec:solution}. Case studies are presented in Section \ref{sec:case_study} followed by the conclusions and limitations in Section \ref{sec:conclusion}.

\section{Problem Definition and System Model}
\label{sec:prob_def_sys_model}
The overall goal of this paper is to optimally design the E-truck fleet composition and associated charging station network. Specifically, a charging service provider decides where to locate new CSs among candidate locations. Additionally, the number of charging ports and substation capacity upgrades (size configurations) are optimized. A fleet operator designs the portfolio of fleet vehicle types, and the optimal routing and charging strategies to deliver customer demands within given time windows while avoiding battery depletion. In this section, we define the problem and present the intuitive illustration to our proposed model.

\subsection{Problem Description}
The problem is defined on a directed graph $\mathcal{G} = \{\mathcal{E}, \mathcal{V}\}$, where $\mathcal{E}$ is the set of all edges\footnote{We will use edge, link, and arc interchangeably.} and $\mathcal{V}$ is the collection of all nodes. Specifically, nodes in $\mathcal{V}$ are categorized into three different types: a depot node $D_0$, customer nodes \{$C_{1},C_{2},...,C_{n}$\} in set $\mathcal{C}$, and candidate charging station (CS) nodes \{$F_{1},F_{2},...,F_{m}$\} in set $\mathcal{F}$. Successive visits of nodes are represented with chosen edges. This is the aforementioned \textit{customer-node based network}. We assume the following common rules: \footnote{Background on the vehicle routing problem and its common formulations can be found in \cite{yang2015battery, hof2017solving, schiffer2016ecvs, schiffer2018designing,schiffer2017adaptive, schiffer2018strategic, schiffer2017electric, erdougan2012green, li2015multiple, hiermann2016electric}.}:
\begin{enumerate}
    \item All customer nodes are visited once and only once by one vehicle during one duty cycle.\label{assump_single_visit}
    \item All E-trucks depart from the depot $D_0$ and return to the same depot after completing the assigned logistic tasks.
    \item Customer demands are represented in the aggregate sense with real values and without specifications, e.g. weight, size, or shape.
\end{enumerate}

\begin{figure}[!tp]
  \begin{center}
  \includegraphics[trim = 10mm 0mm 5mm 0mm, clip, width= 3.5in]{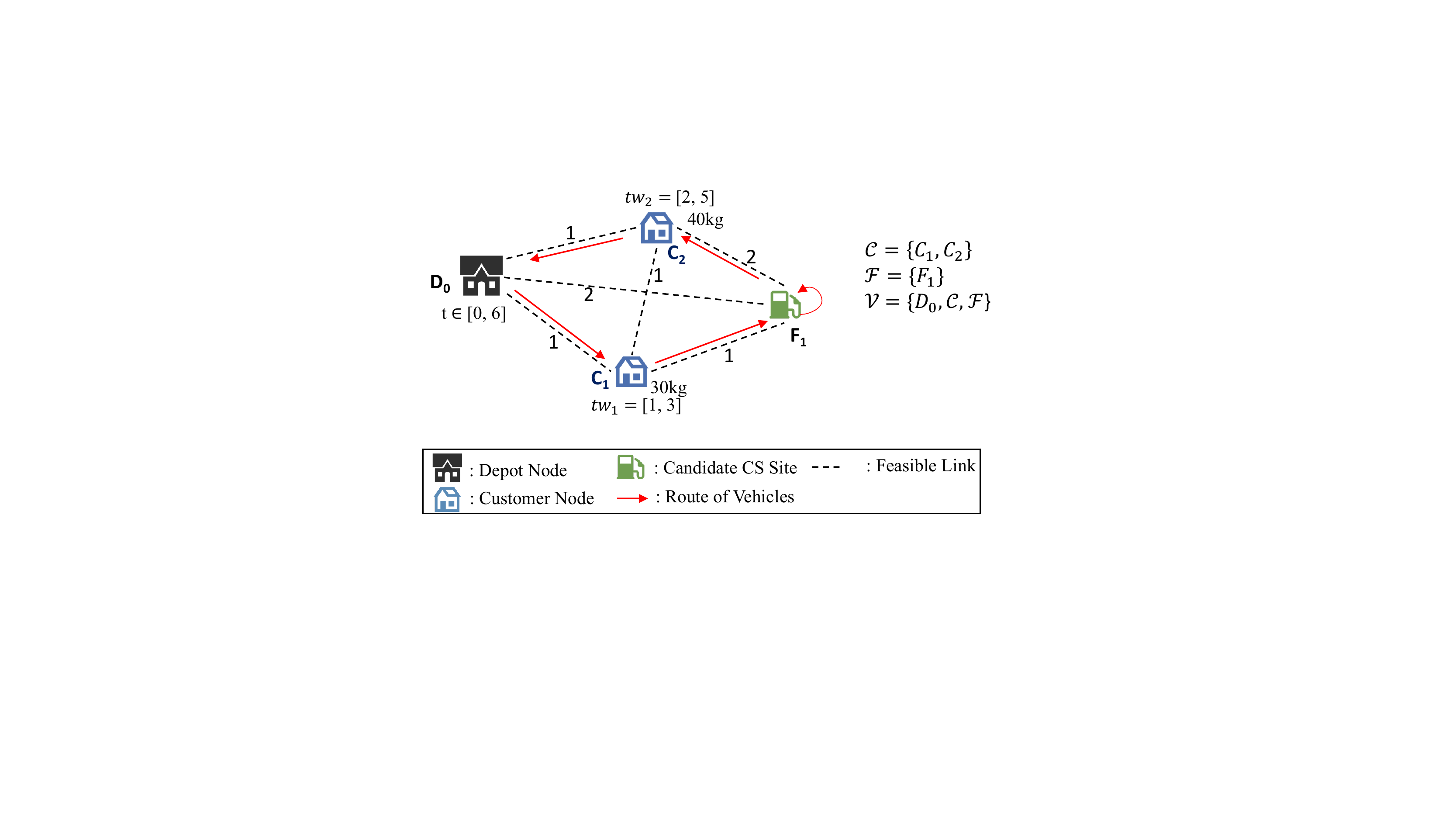}
  \caption{Illustrative network.}\label{fig:net_before_expansion}
  \end{center}
\end{figure}

An illustrative network and toy example is given in Fig.\ref{fig:net_before_expansion}. One depot $D_0$, two customer nodes $C_1$ and $C_2$, and one charging node $F_1$, dashed lines indicate feasible links. Assume there is an E-truck with a driving range of 4 units of length. One possible route is colored in red with arrow directing its trajectory: the E-truck will first make a stop at $C_1$ due to given time window [1,3], then recharge at $F_1$ and go to $C_2$, whose latest required arrival time is 5. Upon task completion, the E-truck will make a return to $D_0$. Alternatively, the feasible routing plan can be chosen as $D_0-C_1-D_0$ and $D_0-C_2-D_0$. Hence, recharging is not required, but two E-trucks are needed to fulfill the task.


In this model, every time index is inherently associated with the node. We loose information to concurrent charging sessions when multiple E-trucks are traversing on the graph. Hence, the model is unable to consider configurations of the charging infrastructure, i.e. the number of ports and transformer upgrades. A work-around is introducing indicator functions to determine specific charging periods; or full state-space-time layered graph is used (Section \ref{subsec:lit_review}). Both approaches either impose nonlinearities or severe scaling issues (number of nodes and links explodes). In the next subsection, we propose a different graph expansion approach to capture the time information neatly.

\subsection{Proposed Model:  Partial Time Expanded Network Model} \label{sec:pten} 
\begin{figure}[!tp]
  \begin{center}
  \includegraphics[scale=0.7]{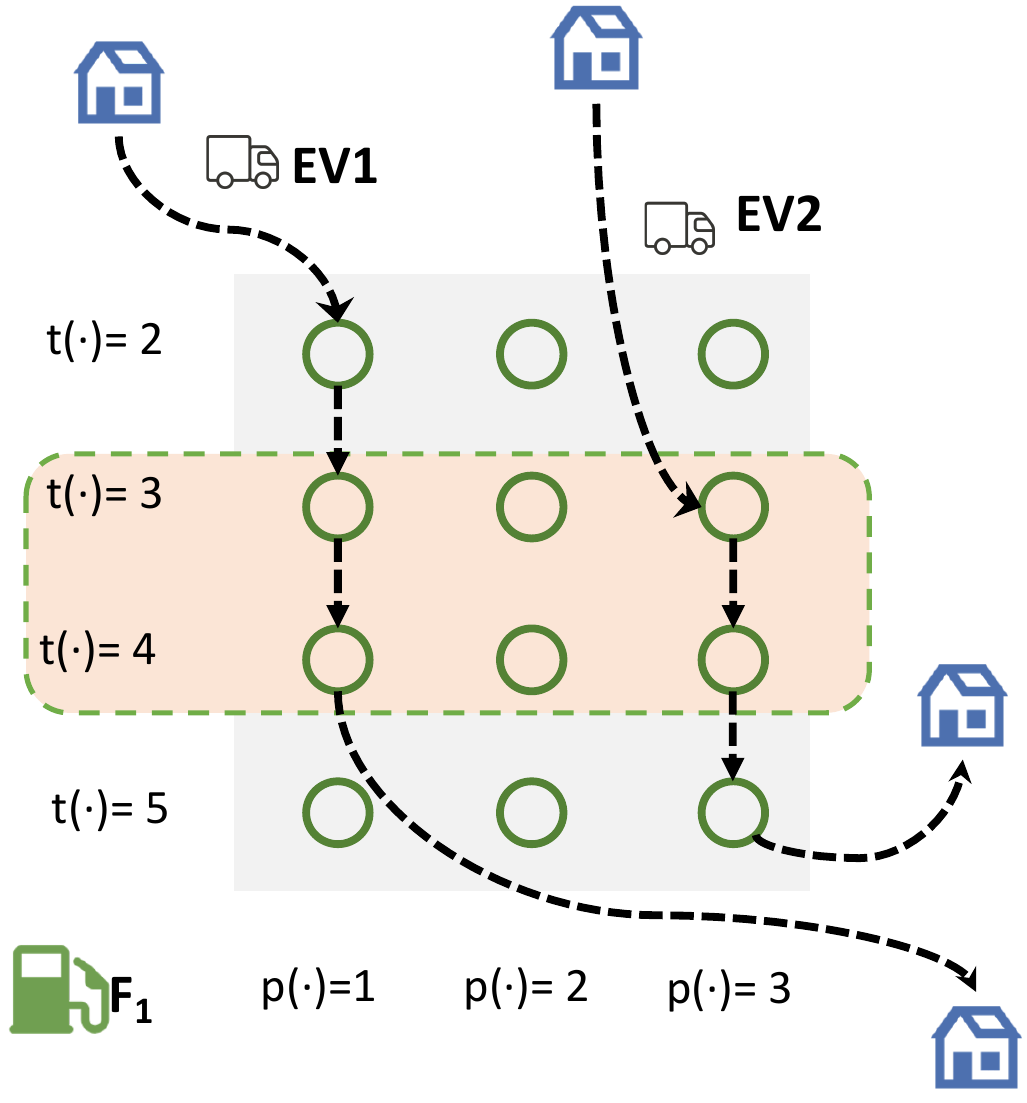}
  \caption{An example for a partial time expanded charging station.}\label{fig:illu_PTEN}
  \end{center}
\end{figure}

We propose to encode the time expansion solely on the charging station nodes, avoiding other unnecessarily added nodes. Namely, this is a partial time expansion. Each original candidate CS node $F_i\in \mathcal{F}$ is expanded across time and charging ports. A two-dimensional time-port graph (Fig.\ref{fig:illu_PTEN}) is introduced to represent a candidate CS node. Altogether, $|\mathcal{T}_i|\cdot|s_{i}|$ dummy nodes\footnote{We will use dummy nodes, dummies, and virtual nodes interchangeably.} are introduced to represent node $F_i$, where $|\mathcal{T}_i|$ is the time horizon (i.e. the number of feasible visiting time slots at $F_i$) and $|s_{i}|$ is the CS size (i.e. the number of charging ports). Each of these nodes, as shown in Fig.\ref{fig:illu_PTEN}, encodes two index trackers: the time availability index $t(\cdot)$ and the charging port index $p(\cdot)$. We denote a set $\mathcal{M}_{i}$ to represent these expanded nodes.

We have defined nodes and now will construct feasible links between nodes within this set $\mathcal{M}_{i}$. A link from node $j$ to node $m$, $\{j,m\in \mathcal{M}_{i}\}$, is created if $t(m)-t(j)=\Delta t$ and $p(m)=p(j)$, indicating an E-truck charges at port $p(j)$ for one time step ($\Delta t$) starting at $t(j)$. We denote the set of all internal links at station node $F_i$ as $\mathcal{A}_{i}$. A subset of time specific links is defined as $\mathcal{A}_{i}(t)= \{(j,m)\in \mathcal{A}_{i} \ | \ t(j)=t, t(m)=t+\Delta t\}$. The original links connecting between customer nodes and the station nodes are reconnected accordingly. With this expansion, real time charging power at station $F_i$ can be easily computed by counting the number of traversed links in $\mathcal{A}_{i}(t)$. Take the case in Fig. \ref{fig:illu_PTEN} as an example, the connection represents that E-truck 1 charges during period $[2,4]$ and E-truck 2 charges during $[3,5]$. Hence, at least two chargers are needed as both E-trucks are present during time $[3,4]$ (box color coded).

\begin{figure}[!tp]
  \begin{center}
  \includegraphics[width=\columnwidth]{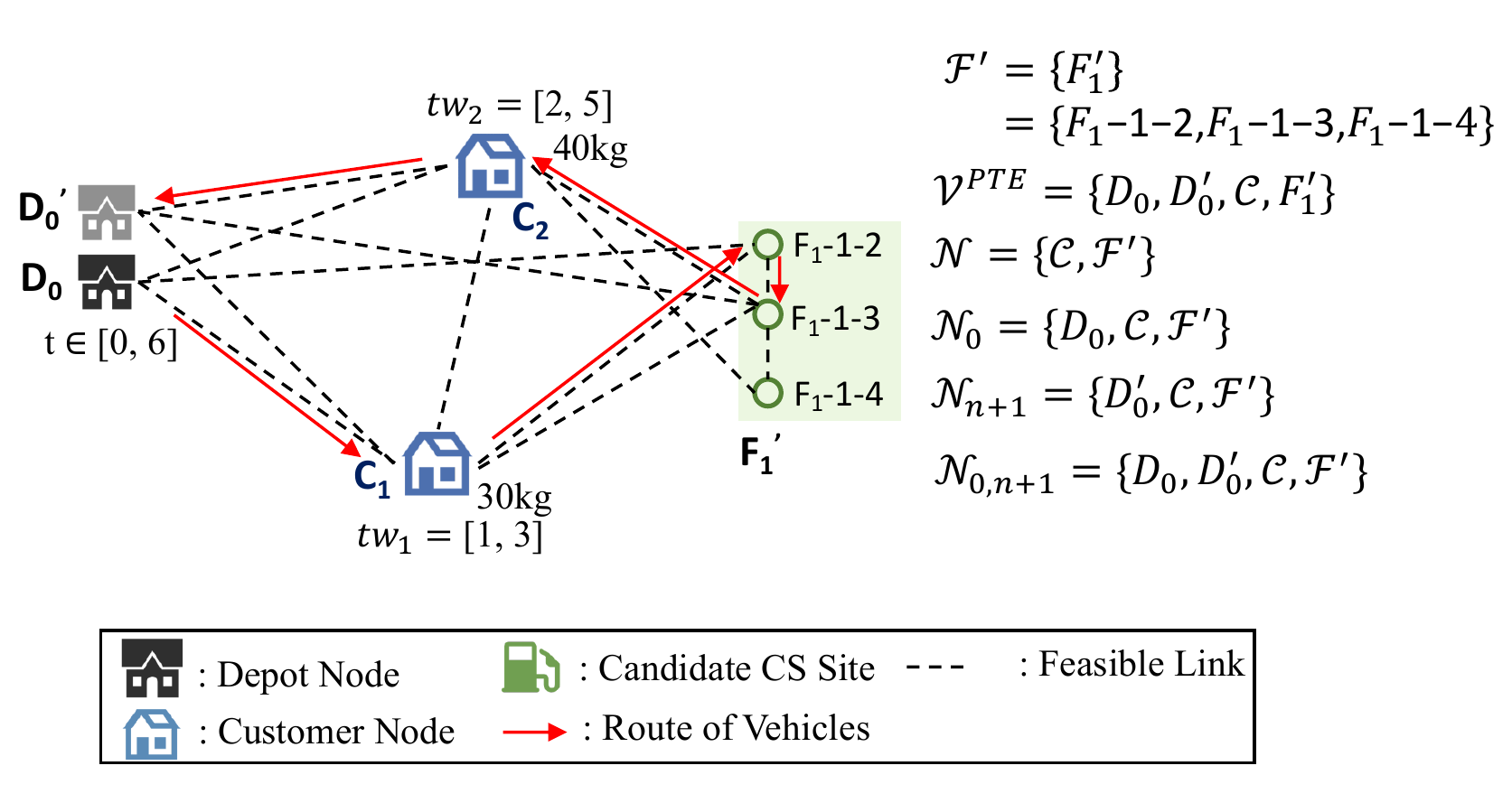}
  \caption{Expanded illustrative network. Here, the notation $F_1-x-y$ represents the $x^{th}$ charger and the corresponding time slot $y$. }\label{fig:net_after_expansion}
  \end{center}
\end{figure}

A copy of the depot node is also created as the sink node $D_0^{\prime}$ (due to assumption \ref{assump_single_visit}, this is a common practice). A corresponding partial time expanded network is presented in Fig.\ref{fig:net_after_expansion}. The charging activity at $F_1$ is then modeled by the link ($F_1-1-2$, $F_1-1-3$). 

We will denote the expanded network as $\mathcal{G}^{PTE} = \{\mathcal{E}^{PTE}, \mathcal{V}^{PTE}\}$. Formal notations are summarized in Table \ref{node_arc_summary}, but relevant sets are also given in Fig.\ref{fig:net_before_expansion} and Fig.\ref{fig:net_after_expansion}. We present the nomenclature in Table \ref{notation} and are now ready to formally introduce the planning problem formulation in the following section. 

\begin{table}[!t]
\renewcommand{\arraystretch}{1.3}
\caption{Nodes and links before and after network expansion}
\label{node_arc_summary}
\centering
\begin{threeparttable}
\begin{tabular}{|c|c|}
\hline
Network & Nodes\\
\hline
$\mathcal{G}$ & $\mathcal{V}=\{\mathcal{D}_0,\mathcal{C},\mathcal{F}\}$ \\
\hline
$\mathcal{G}^{PTE}$ & $\mathcal{V}^{PTE}=\mathcal{V}\cup \mathcal{D}_0^{'}\cup \{\mathcal{F}'=\{\mathcal{M}_i | i\in \mathcal{F}\}\}\setminus \mathcal{F}$\\
\hline
& Arcs\\
\hline
$\mathcal{G}$ &$\mathcal{E}=\{(i,j)| i,j\in\mathcal{V}, i\neq j\}$\\
\hline
\tnote{*}$\mathcal{G}^{PTE}$  &\makecell{$\mathcal{E}^{PTE}=\{(i,j)| i,j\in\mathcal{V}^{PTE}, i\neq j\}$\\$\setminus \{(i,j)| i,j\in\mathcal{F}^{'}_k,t(i)-t(j)\neq \Delta t\}$}\\
\hline
\end{tabular}
\end{threeparttable}
\end{table}

\section{Mathematical Formulation of the Problem as Bi-level Programming} \label{sec:math_model_bi_level}

As mentioned in previous sections, this study aims to capture the interactive dynamics between the charging service provider and the fleet operator. We do not assume cooperation between these two players and hence competition between the two entities is inevitable. 
Today, many transportation logistic companies are considering fleet electrification to reduce total cost of ownership, reduce greenhouse gas emissions, and satisfy upcoming regulations \cite{Sripad2019,Liimatainen2019,Buysse2020}. However, optimally designing the vehicle fleet, routing, and charging strategies remains as open questions. Fleet operators often seek consultation from charging service providers. This naturally leads to a leader-follower setting, in which the charging service provider is the leader\footnote{The leader knows the cost function mapping of the follower in this game. The follower, on the other hand, observes the strategies from the leader and always has to take them into account.}. Leader-follower games are also known as Stackelberg games \cite{cruz1975survey}. 

The decision variables and notation for the problem are summarized in Table \ref{notation}. Next, we detail the model for each player as well as for the complete problem.
\begin{table}[!t]
\caption{Notation summary (Alphabetical order)}
\label{notation}
\centering
\begin{tabular}{|p{0.1\linewidth} | p{0.8\linewidth}|}
\hline
\multicolumn{2}{|c|}{\textbf{Variables}}\\
\hline
$b_i^k$& Continuous variable, battery level (kWh) of type $k$ vehicle when arrived at node $i$\\
\hline
$\Delta P_i$ & Continuous variable, upgrade capacity of the transformer connecting node $i$\\
\hline
$q_i^k$&Continuous variable, load level (kg) of type $k$ vehicle when arrived at node $i$\\
\hline
$s_i$ & Integer variable, number of charging ports installed in the charging station at node $i$\\
\hline
$\tau_{i}$& Integer variable, time when a vehicle arrived at node $i$\\
\hline
$x_{ij}^k$&Binary variable, whether vehicle of type $k$ visits node $j$ after node $i$\\
\hline
$y_i$ & Binary variable, whether a charging station is constructed at node $i$\\
\hline
\multicolumn{2}{|c|}{\textbf{Sets}}\\
\hline
$\mathcal{E}^{N}$& Sets of links in the expanded network excluding the internal links within any CS, i.e. $\mathcal{E}^{PTE}\setminus(\cup \mathcal{A}_{i})$\\
\hline
$\mathcal{F}'$&Set of all dummy CS nodes, i.e. $\{\mathcal{M}_i | i\in \mathcal{F}\}$\\
\hline
$\mathcal{K}$& Set of vehicle type index\\
\hline
$\mathcal{M}_{i}(t) $& Set of dummy nodes for CS $i$ with whose time tracker equals $t$, i.e. $\{j|j\in\mathcal{M}_{i}, t(j)=t\}$\\
\hline
$\mathcal{N} $& Set of all dummy CS nodes and customer nodes, i.e. $C \cup \mathcal{F}'$\\
\hline
$\mathcal{N}_0 $& Set of all nodes except the sink depot, i.e. $\mathcal{N} \cup \{D_0\}$\\
\hline
$\mathcal{N}_{n+1} $& Set of all nodes except the source depot, i.e. $\mathcal{N} \cup \{D_0^{\prime}\}$\\
\hline
$\mathcal{N}_{0,n+1} $& Set of all nodes, i.e. $\mathcal{N}\cup \{D_0, D_0^{\prime}\}$\\
\hline
\multicolumn{2}{|c|}{\textbf{Parameters}}\\
\hline
$B^k$ &  Battery capacity of a type $k$ truck, in [kWh]\\
\hline
$c_i^p$ & Cost of one charging port at node $i$, in [\$/port]\\
\hline
$c_i^s$ & Cost of substation capacity upgrade at node $i$, in [\$/kW]\\
\hline
$c^{e}_i$ & Service fee charged by the charging station $i$, in [\$/kWh]\\
\hline
$c^{k}$ & Cost of vehicle type $k$, in [\$]\\
\hline
$c_{ij}^{k}$ & Cost of travel from node $i$ to node $j$ of type $k$, in [\$]\\
\hline
$d_i$ & Customer demand at node $i$\\
\hline
$d_{ij}$ & Distance between node $i$ and node $j$, in [km]\\
\hline
$\Delta t$ &  Time step\\
\hline
$\pi_{i,t}$ & Electricity price at node $i$ at time $t$, in [\$/kWh]\\
\hline
$p^\text{rated}$ & Charger rated power, in [kW]\\
\hline
$P_{i,t}$ & Substation availability at time $t$, at node $i$, in [kW]\\
\hline
$Q^k$ &  Freight capacity of a type $k$ truck, in [kg]\\
\hline
$r^k$ &  Energy consumption rate of type $k$, in [kWh/km]\\
\hline
$t_i^e$, $t_i^l$ & Earliest arrival/Latest departure time at node $i$\\
\hline
$t_i^s$ & Required service time of node $i$\\
\hline
$t_{ij}$ & Travel time from node $i$ to node $j$\\
\hline
$\zeta_{s/v}$ & Capital recovery factors for the charging station and the vehicles respectively\\
\hline
\end{tabular}
\end{table}

\subsection{CSP's Problem: Charging Network Design and Operation}
The leader CSP aims to minimize its overall costs by optimally placing and sizing the new charging stations. The binary variable $y_i$ is used to indicate the construction decision at the specific site $F_i$, and the integer variable $s_i$ represents the number of chargers to be installed at $F_i$.

The overall cost $g^{L}$ constitutes two parts. The first part is the capital expenditure (CAPEX), namely the costs for installing ports and for upgrading the local transformer if necessary. The second part is the operational profit introduced by providing charging service to E-trucks with the predetermined service fee $c_i^e$ (\$/kWh). Mathematically, $g^{L}$ is expressed as
\begin{equation}\label{CSP_obj}
\begin{split}
   g^{L}(\y, \s, \DeltaP) & = \sum_{i \in \mathcal{F}} \{\zeta_{s} \cdot (c_{i}^{s}s_{i}+c_{i}^{p}\Delta P_i) - c^{e}_i E^{ch}_{i}\}.
   \end{split}
\end{equation}




Here, $E^{ch}_{i}$ represents the total electricity delivered to E-trucks at station $F_i$. Based on the proposed PTEN model in Section \ref{sec:pten}, it is calculated as 
\begin{equation}
    E^{ch}_{i} = p^\text{rated} \Delta t\cdot \sum_k \sum_{(j,m) \in\mathcal{A}_i} x_{j,m}^{k}.
\end{equation}
Factor $\zeta_{s}$ in $g^{L}$ converts the life-cycle fixed cost into its annual equivalent level, which is calculated as 
\begin{equation}
	\zeta_{s} = \frac{r(1+r)^{Y_{s}}}{(1+r)^{Y_{s}}-1},
\end{equation}
where $r$ is the cash discount rate and $Y_{s}$ is the service life of the charging station.

When building the charging stations, the CSP should ensure adequate chargers and sufficient transformer capacity to supply power to the visited vehicles, which leads to
\begin{equation}\label{charging_cons}
    \sum_{k\in \mathcal{K}}\sum_{(j,m)\in\mathcal{A}_{i}(t)} x_{j,m}^{k} \leq s_i \quad \forall i \in \mathcal{F}, \forall t \in \mathcal{T}_i,
\end{equation}
\begin{equation}\label{upgrade_cons}
	s_{i} \cdot p^\text{rated} - y_{i} \cdot P_{i,t} \leq \Delta P_i   \quad \forall i \in \mathcal{F}.
\end{equation}
There are also constraints on the station size and variable domain constraints: 
\begin{equation}\label{size_cons}
    y_{i}s_{i}^\text{min}\leq s_{i}\leq y_{i}s_{i}^\text{max} \quad \forall i \in \mathcal{F},
\end{equation}
\begin{equation}\label{y_cons}
    y_{i} \in\{0,1\}, s_{i} \in \mathbf{Z}^{+}, \Delta P_i\in \mathbf{R}^{+},  \quad \forall i \in \mathcal{F}.
\end{equation}


\subsection{FO's problem: Fleet Design and Operation}
The FO's goal is to decide its E-truck fleet composition and routing plans so that its overall cost is minimized. The cost objective for the FO is
\begin{align}\label{FO_obj}
	g^{F}&(\x)= \zeta_{v}\sum_{k\in \mathcal{K}} \sum_{j \in \mathcal{N}} c^{k} x_{D_{0} j}^{k}+ \sum_{k\in \mathcal{K}} \sum_{i \in \mathcal{N}_{0}} \sum_{j \in \mathcal{N}_{n+1}} c_{i j}^{k} x_{i j}^{k} \nonumber \\
    &+\sum_{i \in \mathcal{F}}\sum_{t \in \mathcal{T}_i} (\pi_{i,t}+c^{e}_i)\sum_k \sum_{(j,m) \in\mathcal{A}_{i}(t)} p^\text{rated} \Delta t \cdot x_{j,m}^{k}.
\end{align}
The first term in $g^{F}$ represents the total E-truck purchase cost, which is converted into the equivalent annual level using $\zeta_{v}$. The second term yields the traveling cost of the fleet. The last term calculates the cost of charging, where the per unit charging cost involves the electricity price $\pi_{i,t}$ plus the service fee $c^{e}_i$ posed by the CSP.

Vehicle routing must respect resource constraints along the network, including time windows, payload capacity, energy, etc., which are given as follows.

\subsubsection{Network flow constraints}
\begin{equation}\label{single_visit_cust}
	\sum_{k\in \mathcal{K}} \sum_{j \in N_{n+1}, j \neq i} x_{i j}^{k}=1 \quad \forall i \in \mathcal{C}, 
\end{equation}
\begin{equation}\label{visit_cs}
	\sum_{k\in \mathcal{K}} \sum_{j \in N_{n+1}, j \neq i} x_{i j}^{k} \leq 1 \quad \forall i \in \mathcal{F}^{\prime},
\end{equation}
\begin{equation}\label{flow_conser}
    \sum_{j \in N_{0}, j \neq i} x_{j i}^{k}-\sum_{j \in N_{n+1}, j \neq i} x_{i j}^{k}=0 \quad \forall i \in N, \forall k \in \mathcal{K},
\end{equation}

Constraint (\ref{single_visit_cust}) requires each customer to be visited once and only once, while for the expanded charging station nodes, this requirement is relaxed in (\ref{visit_cs}). Flow conservation of each node, except the source and sink, is expressed by (\ref{flow_conser}).

\subsubsection{Time window constraints}
\begin{equation}\label{time_vist}
	t_{i}^{e} \leq \tau_{i} \leq t_{i}^{l} \quad \forall i \in \mathcal{C}\cup D_{0}, D_{0}^{'},
\end{equation}
\begin{equation}\label{time_to_cust}
\begin{split}
	\tau_{j}-\tau_{i} \geq\left(t_{i j}+t_{i}^{s}\right) x_{i j}^{k}-\left(1-x_{i j}^{k}\right) T \\
	\quad \forall i \in \mathcal{N}_{0}, \forall j \in \mathcal{N}_{n+1}, \forall k \in \mathcal{K}, j \neq i, 
\end{split}
\end{equation}
\begin{equation}\label{time_to_cs1}
\begin{split}
	t(j)-\tau_{i} \geq\left(t_{i j}+t_{i}^{s}\right) x_{i j}^{k}-\left(1-x_{i j}^{k}\right) T\\
	\quad \forall i \in \mathcal{N}_{0}, \forall j \in \mathcal{F}^{\prime}, \forall k \in \mathcal{K}, j \neq i,
\end{split}
\end{equation}
\begin{equation}\label{time_to_cs2}
\begin{split}
	t(j)-\tau_{i} \leq\left(t_{i j}+t_{i}^{s}\right) x_{i j}^{k}+\left(1-x_{i j}^{k}\right) T \\
	\quad \forall i \in \mathcal{N}_{0}, \forall j \in \mathcal{F}^{\prime}, \forall k \in \mathcal{K}, j \neq i, 
\end{split}
\end{equation}
\begin{equation}\label{time_to_cs3}
	\tau_j = \sum_{k\in \mathcal{K}}\sum_{i} x_{i j}^{k} \cdot t(j) \quad \forall i \in \mathcal{N}_{0}, \forall j \in \mathcal{F}^{\prime}
\end{equation}

The arrival time of a vehicle at a customer point $i$ must respect the customer's service time window $\left[t_{i}^{e},t_{i}^{l}\right]$ \eqref{time_vist}. For the depots, the time window is set as $\left[0,T\right]$. Constraint \eqref{time_to_cust} expresses the relationship between two successive customer nodes $i$ and $j$ at their respective visited times. When $x_{ij}^{k}=1$, then E-truck $k$'s arrival time at customer $j$ depends on the traveling time between $i,j$ and the service time at $i$. However, when $x_{ij}^{k}=0$, i.e. $j$ is not visited after $i$, then this constraint is relaxed.

Constraints (\ref{time_to_cs1})-(\ref{time_to_cs3}) describe the evolution of visiting times when an E-truck is driving towards a charging station node. Since each dummy node is strictly associated with one specific time slot, the corresponding relation $x_{ij}^{k} = 1$ is true only if the arrival time $\tau_{j}$ at the charging node $j$ matches $t(j)$, as shown in (\ref{time_to_cs1}). Again,  those constraints are relaxed if $x_{ij}^{k}=0$.

\subsubsection{Freight capacity constraints}
\begin{align}
    \label{freight_cap}& q_{j}^{k} \leq q_{i}^{k}-d_{i} x_{i j}^{k}+\left(1-x_{i j}^{k}\right) Q^{k}\\
    \nonumber &\quad \forall i \in \mathcal{N}_{0}, \forall k \in \mathcal{K}, \forall j \in \mathcal{N}_{n+1}, j \neq i,\\
    \label{freight_level}& 0 \leq q_{i}^{k} \leq Q^{k} \quad \forall k \in \mathcal{K}, \forall i \in \mathcal{N}_{0, n+1},
\end{align}

Following the same modeling philosophy from above, the available freight loads at each node along the route are tracked using \eqref{freight_cap}. Constraint \eqref{freight_level} ensures that the E-trucks are never overloaded.

\subsubsection{Energy consumption/recharge constraints}
\begin{equation}\label{batt_consump}
\begin{split}
	b_{j}^{k} \leq b_{i}^{k}-r^{k} d_{i j} x_{i j}^{k}+\left(1-x_{i j}^{k}\right) B^{k} \\
	\quad \forall (i,j) \in \mathcal{E}^{N}, \forall k \in \mathcal{K}, j \neq i
\end{split}
\end{equation}
\begin{equation}\label{batt_recharge}
\begin{split}
	b_{j}^{k} \leq b_{i}^{k}+p^\text{rated} \cdot \Delta t\cdot x_{i j}^{k}+\left(1-x_{i j}^{k}\right) B^{k} \\
	\quad \forall (i,j) \in \mathcal{A}_z, \forall z\in \mathcal{F}, \forall k \in \mathcal{K}
\end{split}
\end{equation}
Given the limited range of E-trucks, it is crucial to track available battery energy while traveling, which is modeled by \eqref{batt_consump} and \eqref{batt_recharge}. These constraints are relaxed when $x_{ij}^k=0$ by using the term $(1-x_{i j}^{k})B^k$. We assume a constant energy consumption rate while traversing to customers and a constant charging rate while traversing the charging links. Visiting consecutive charging nodes at one physical location represents charging for multiple time slots.


\begin{equation}\label{batt_at_0}
	b_{0}^{k}=B^{k} \quad \forall k \in \mathcal{K}, 
\end{equation}
\begin{equation}\label{batt_level}
	0 \leq b_{i}^{k} \leq B^{k} \quad \forall k \in \mathcal{K}, \forall i \in \mathcal{N}_{n+1},
\end{equation}
\begin{equation}\label{batt_charge}
	b_{i}^{k}+p^\text{rated} \cdot \Delta t \leq B^{k} \quad \forall k \in \mathcal{K}, \forall i \in \mathcal{F}^{\prime},
\end{equation}
We assume all E-trucks start fully charged at depot \eqref{batt_at_0}. The battery is never depleted nor overcharged as enforced by \eqref{batt_level} and \eqref{batt_charge}, respectively.

\subsubsection{Simultaneous charging constraint}
\begin{align}\label{couple_cons}
     \sum_{k\in \mathcal{K}}\sum_{(j,m)\in\mathcal{A}_{i}(t)} x_{j,m}^{k} \leq s_i \quad \forall i \in \mathcal{F}, \forall t \in \mathcal{T}_i.
\end{align}
The number of simultaneous charging E-trucks must respect the physical charging station size limit. 

\subsubsection{Additional variable domains}
\begin{align}\label{var_dom_fo}
    & x_{i,j}^{k} \in\{0,1\} \quad \forall i \in \mathcal{N}_0,  \forall j \in \mathcal{N}_{n+1}\\
    &\nonumber \tau_{i}\in \mathbf{Z}^{+},\quad q_{i}, b_{i} \in \mathbf{R}^{+}\quad \forall i \in \mathcal{N}_{0,n+1}.
\end{align}

\subsection{Joint Problem as a Stackelberg Game}
Given the CSP (leader) and FO (follower) optimization models above, we now integrate them to yield the complete joint planning problem: 
\begin{subequations}
\label{complete_problem}
\begin{flalign}
    \min_{\y, \s, \DeltaP} \quad & g^{L}(\y,\s,\DeltaP;\boldsymbol{x^{*}})\\ 
     \mbox{s. to:}   \quad & h^{L}(\y,\s,\DeltaP)\leq 0\\
            \quad & (\boldsymbol{x^{*},\tau^{*},b^{*},q^{*}})=\mathop{\arg\min}_{\x,\rtau,\rb,\q} \quad g^{F}(\x)\label{lower_level_problem}\\
            & \nonumber\qquad \qquad \mbox{s. to:} \quad h^{F}(\s, \x, \rtau, \rb, \q)\leq 0.
\end{flalign}
\end{subequations}
The constraint set $h^{L}(\y,\s,\DeltaP)$ contains \eqref{charging_cons}-\eqref{y_cons} and $h^{F}(\s, \x, \rtau, \rb, \q)$ includes \eqref{single_visit_cust}-\eqref{var_dom_fo}. For the reader's convenience, we have colored the leader's optimization variables blue and the follower's optimization variables red. Black bold variables are fixed optimization variables \footnote{We would emphasize that although this is a fleet sizing, facility siting and sizing, and vehicle routing joint decisions, the model is also able to consider the current existing charging network. We simply convert the corresponding decision variables to input parameters. It is a degenerate case of our model. }.

A key benefit of the proposed PTEN is the overall mathematical formulation \eqref{complete_problem} maintains a mixed integer linear programming structure. However, solving this model is still highly non-trivial. The main challenges are twofold: (\textrm{i}) The overall model is a bi-level mixed integer problem (Bi-MILP) and integer variables exist in both the upper and lower levels. In this case, the commonly-used KKT (Karush-Kuhn-Tucker)-based single-level reformulation method is not applicable. (\textrm{ii}) The electric vehicle routing problem, embedded as the essential part of the overall problem, is an NP-hard problem whose scale grows dramatically with the size of the network. This holds true even with the partial time expansion, which mitigates but does not eliminate the computational complexity. Additionally, off-the-shelf solvers cannot be directly applied, even for a medium sized system. 

To address the above challenges, a customized double loop solution algorithm is proposed and presented in the next section.

\section{Solution Algorithm Design} \label{sec:solution}
\begin{figure*}[!tp]
  \begin{center}
  \includegraphics[scale=0.7]{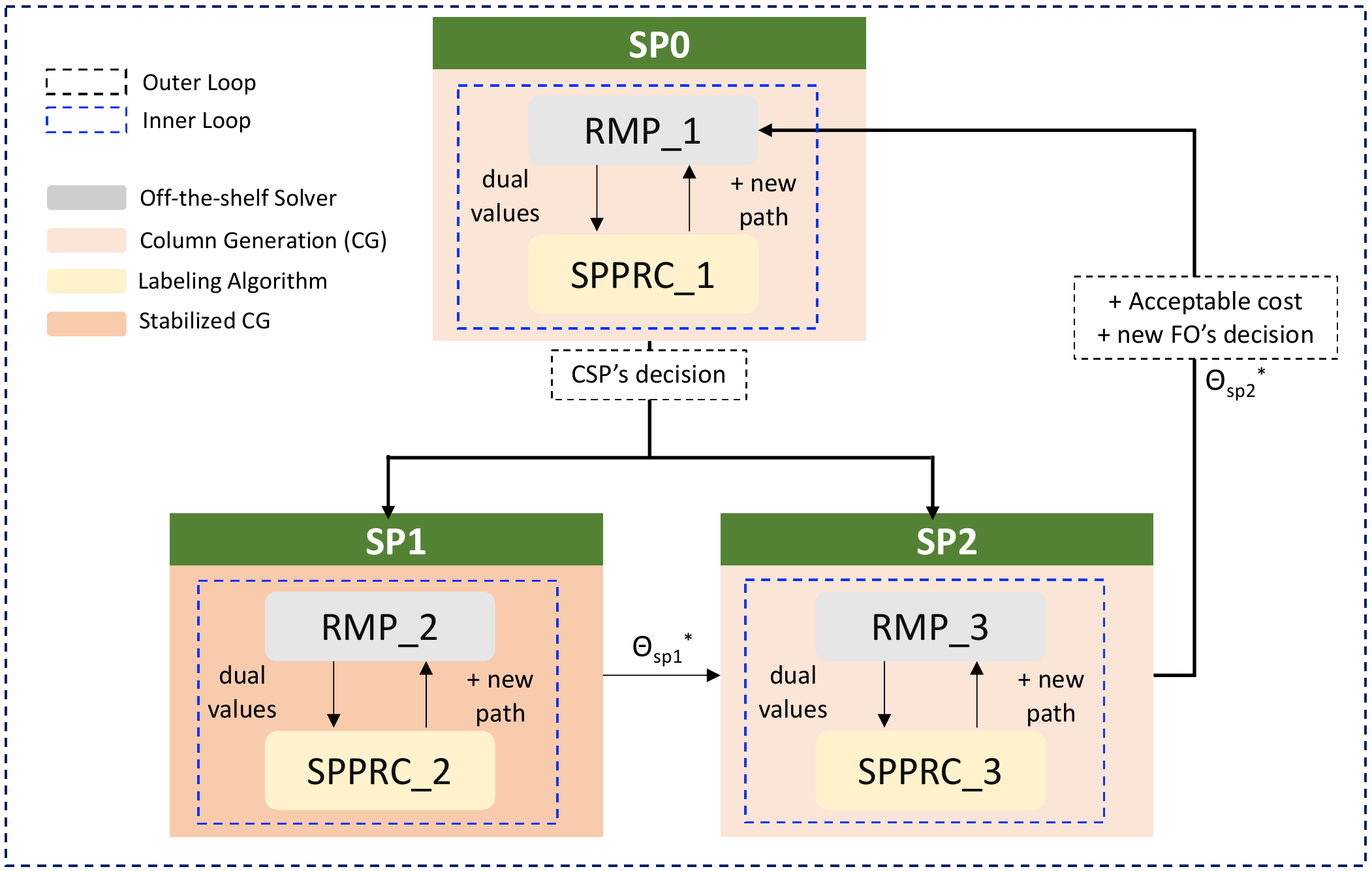}
  \caption{The overall solution architecture. Notice that (SP2) has been detailed in Section \ref{sec:bi_level_decomposition}. It serves as a CSP feasibility check to evaluate the solution set from its successor (SP1). Its reformulation follows the exact same process from Section \ref{sec:column_generation}.RMP and SPPRC denote the restricted mater problem and the shortest path problem with resource constraints respectively, which are detailed in the inner-loop design in Section \ref{sec:column_generation}}
  \label{overall_algo}
  \end{center}
\end{figure*}

As the network size increases, model \eqref{complete_problem} becomes a massive Bi-MILP problem that the off-the-shell solver is not capable to handle. Here, we propose a double-loop framework to release computation burden and iteratively resolve this issue. To offer a more straightforward visualization of the overall architecture, we present a diagram flow in Fig.\ref{overall_algo}.

\subsection{Outer-loop: Solving Bi-MILP with reformulation and decomposition} \label{sec:outer-loop}
To address the first challenge (i), we facilitate the overall computation of the Bi-MILP with a reformulation and decomposition method \cite{zeng2014solving}. We denote this as an outer loop design. The main idea of this approach and the problem-specific implementation are presented below and readers are referred to the original references for further theoretical details.

\subsubsection{Reformulation}
A key observation is that when all discrete variables in the lower-level FO problem are fixed, then the FO \eqref{lower_level_problem} becomes a pure linear program with continuous decision variables. Then the optimal solution can be represented using KKT conditions. Given a specific realization for the $z$-indexed combination of the lower level discrete variables $\x^z, \rtau^z$, we denote the corresponding KKT conditions as $\Lambda(\x^z, \rtau^z)$. Therefore, if one enumerates all possible combinations of the follower's discrete decisions, we can denote their collection using index set $\mathcal{L}_\text{full}=\{1, \dots, l^\text{max}\}$ where $z \in \mathcal{L}_\text{full}$. Then the original problem can be equivalently formulated as a single-level problem:
\begin{align}
   \textbf{(P0)}\quad  \min & \quad g^{L}(\y, \s, \DeltaP, \x^0) \label{p0_obj}\\
   \mbox{s. to:}\quad & h^{L}(\y, \s, \DeltaP)\leq 0,\\
   &  h^{F}(\s, \x^0, \rtau^0, \rb^0, \q^0)\leq 0, \\
   & \forall z\in \mathcal{L}_\text{full} \nonumber \\ 
   & g^{F}(\x^0)\leq g^{F}(\boldsymbol{x^{z}}), \label{better_obj}\\
   & \y,\s,\DeltaP, \boldsymbol{b^{z},q^{z}} \in \Lambda(\boldsymbol{x^{z},\tau^{z}}). \label{KKT_set}
\end{align}
Variables $ \x^{0}, \rtau^{0},\rb^{0},\q^{0}$ are duplications of the follower's decisions. Constraint (\ref{better_obj}) requires that the FO's objective is at least the same, if not improved, from the discrete (and previous as we shall see) realization $\boldsymbol{x^z}$. Notably, in our case the discrete variables $\boldsymbol{x,\tau}$ uniquely define the routes of the fleet. Once all routes are realized, then values of $\boldsymbol{b,q}$ are implicitly determined. The complexity of (\ref{KKT_set}) can thus be largely reduced.
\subsubsection{Decomposition} \label{sec:bi_level_decomposition}
Instead of directly solving the complete problem \textbf{(P0)} with all possible combinations of $\{\boldsymbol{x,\tau,b}\}$ enumerated, one may solve the problem with a subset $\mathcal{L}^\text{sub}$ of these combinations, i.e. $z \in \mathcal{L}_\text{sub} \subseteq \mathcal{L}_\text{full}$ and gradually enlarge the set. As explained in \cite{zeng2014solving}, the solution of \textbf{(P0)} can be obtained by iteratively solving the following decomposed parts,
\begin{itemize}
    \item A restricted version of \textbf{(P0)} with the subset $\mathcal{L}_\text{sub} \subseteq \mathcal{L}_\text{full}$, denoted as \textbf{(SP0)}. Since only a subset of all the constraints are considered, \textbf{(SP0)} provides a lower bound to the original problem \textbf{(P0)}.
    \item Subproblem 1 \textbf{(SP1)} finds the follower's corresponding best response \{$\boldsymbol{x^{*},\tau^{*},b^{*},q^{*}}$\} to the leader's decisions $\{\boldsymbol{y^{*},s,^{*}\Delta P^{*}}\}$ from \textbf{(SP0)}, i.e.
    \begin{align}\label{SP1}
    {\textbf{(SP1)}} \quad 
    & \mathop{\min}_{\x, \rtau, \rb, \q} \quad g^{F}(\x)\\ 
    \mbox{s. to:}\quad & h^{F}(\boldsymbol{s}^*, \x, \rtau, \rb, \q)\leq 0,
    \end{align}
    \item Subproblem 2 \textbf{(SP2)} performs a feasibility check\footnote{Sometimes there may be multiple non-unique lower-level optimal solutions given the upper-level decision. By solving \textbf{(SP2)}, we select the follower solution that is most in favor of the leader.} and is defined as 
    \begin{align}\label{SP2}
    {\textbf{(SP2)}} \quad \nonumber \min_{\x, \rtau, \rb, \q} & \quad g^{L}(\boldsymbol{y^{*},s^{*},\Delta P^{*}}; \x)\\ 
    \mbox{s. to:}\quad & g^{F}(\x)\leq \theta^{*}_\text{sp1},\\
    & \nonumber h^{F}(\boldsymbol{s^{*}},\x,\rtau,\rb,\q)\leq 0,
    \end{align}
    where $\theta^{*}_{sp1}$ is the optimal value from \textbf{(SP1)}. When a solution is found feasible in \textbf{(SP2)}, then the decision set \{$\boldsymbol{\tilde{x}^{*},\tilde{\tau}^{*},\tilde{b}^{*},\tilde{q}^{*}}$\} represents the most favorable follower action for the leader. We then add it into $\mathcal{L}_\text{sub}$ in \textbf{(SP0)} for the next round of iteration. \textbf{(SP2)} provides $\theta^{*}_{sp2}$ as an upper bound for \textbf{(P0)}, since it clearly finds a feasible solution. 
\end{itemize}

The pseudo code is detailed in Algorithm \ref{algo:bi_milp}:
 \begin{algorithm}[ht]
 \caption{Column and Constraint Generation Algorithm for the Joint Planning Bi-MILP}
 \begin{algorithmic}[1]
 \renewcommand{\algorithmicrequire}{\textbf{Input:}}
 \renewcommand{\algorithmicensure}{\textbf{Output:}}
 \REQUIRE model parameters and convergence margin $\epsilon$
 \ENSURE  optimal solution for both CSP and FO
 \\ \textit{- Initialization} :
  \STATE  Set $LB=-\infty$, $UB=\infty$, and $l=0$
 \\ \textit{- Loop Process :}
  \WHILE {$UB-LB > \epsilon$}
  \STATE Solve (SP0) with current combinations $z=1, \dots, l^\text{sub}$ and obtain \{$\boldsymbol{y^{*}_z,s^{*}_z,\Delta P^{*}_z}$\} and the optimal objective $\Theta^{*}_\text{RMP}$, set $LB=\Theta^{*}_\text{RMP}$.
  \STATE Solve (SP1) given \{$\boldsymbol{y^{*}_z,s^{*}_z,\Delta P^{*}_z}$\} as fixed, and obtain \{$\boldsymbol{x^{*},\tau^{*},b^{*},q^{*}}$\} and the optimal objective as $\theta^{*}_\text{sp1}$
  \STATE Solve (SP2) with \{$\boldsymbol{y^{*}_z,s^{*}_z,\Delta P^{*}_z}$\} as fixed, and obtain \{$\boldsymbol{\tilde{x}^{*},\tilde{\tau}^{*},\tilde{b}^{*},\tilde{q}^{*}}$\} and the optimal objective as $\theta^{*}_\text{sp2}$
  \IF {(SP2) is feasible}
  \STATE  Set \{$\boldsymbol{x^{z+1},\tau^{z+1}}$\} as \{$\boldsymbol{\tilde{x}^{*},\tilde{\tau}^{*}}$\},\\ $UB = \min \{UB,\theta^{*}_\text{sp2}$\} 
  \ELSE 
  \STATE Set \{$\boldsymbol{x^{z+1},\tau^{z+1}}$\} as \{$\boldsymbol{x^{*},\tau^{*}}$\}
  \ENDIF
  \STATE Add the new optimal cut corresponding to \{$\boldsymbol{x^{z+1},\tau^{z+1}}$\} to the (SP0), set $z=z+1$.
  \ENDWHILE
 \RETURN ${y^{*}_z,s^{*}_z,\Delta P^{*}_z, x^{z},\tau^{z},b^{z}}$
 \end{algorithmic}
 \label{algo:bi_milp}
 \end{algorithm}
While challenge (i) regarding integer variables at both the upper and lower levels has been addressed, challenge (ii) remains. This motivates the inner loop algorithm design in the next subsection.

\subsection{Inner loop: Solving E-V/LRP with column generation and labeling algorithms} \label{sec:column_generation}
Notice that the three subproblems are structurally similar and the electric vehicle routing problem serves as the core in \textbf{(SP1)} and \textbf{(SP2)} and the location planning is encoded in  \textbf{(SP0)}. To solve these problems efficiently, we first reformulate them as generalized set-partitioning problems and again solve in an iterative fashion to reduce computation burden. 
\subsubsection{Set partitioning reformulation}
We start with a reformulation of \textbf{(SP1)}, which only involves the E-VRP. It is given by
\begin{align}
    \textbf{(MP1)} \quad & \min_{\iota^k} \sum_{k \in \mathcal{K}} \sum_{r \in \Omega_k} c_r^k \cdot \iota^k \label{obj:label_mp}\\
    \mbox{s. to:} \quad & \sum_{k \in \mathcal{K}} \sum_{r \in \Omega_k} a_{i r}^k \iota^k=1, \quad \forall i\in\mathcal{C}\quad :\boldsymbol{\gamma} \label{constr:mp_route1}\\
    & \iota^k \in\{0,1\}, \forall r \in \Omega_{k}, \forall k \in \mathcal{K} \label{constr:mp_domain}
\end{align}
where $\Omega_k$ is the collection of all feasible routes for vehicles of type $k$ satisfying constraints (\ref{single_visit_cust})-(\ref{var_dom_fo}). Parameter $c_r^k$ is the corresponding cost of route $r$, including vehicle cost, traveling cost and charging cost. A binary variable $\iota^k$ is introduced as the route selection indicator. Parameter $a_{i r}^k$ indicates whether a customer node $i$ is visited within the current route. The objective (\ref{obj:label_mp}) is thus to find a set of selected routes that minimizes the total cost subject to the customer visiting constraint (\ref{constr:mp_route1}) and the variable domain constraints (\ref{constr:mp_domain}).

Compared with the pure E-VRP problem in \textbf{(SP1)}, the reformulated subproblem for the E-LRP problem in \textbf{(SP0)} is more complex as it involves decision variables for charging station siting and sizing. Nevertheless, the key technique is also to re-express the formulation as a function of routes. We introduce two new parameters $b_{j r}^k$ and $co_r^k$ to indicate the visits of charging nodes and the CSP's charging service revenue (negative cost) along the route $\iota^k$, respectively. The reformulated model \textbf{(MP0)} for the leader CSP is:
\begin{subequations}
\label{obj:label_mp-lrp}
\begin{flalign}
\textbf{(MP0)} \quad & \min_{\iota^k,s_i,\Delta P_i} \sum_{i \in \mathcal{F}} \zeta_{s} \cdot (c_{i}^{s}s_{i}+c_{i}^{p}\Delta P_i) \label{obj_mp0-invest}\\
    &\qquad +\sum_{k \in \mathcal{K}} \sum_{r \in \Omega_k} co_r^k \cdot \iota^k \label{obj_mp0-energy}
    \end{flalign}
\end{subequations}
\begin{align}
    \mbox{s. to:} \quad & \sum_{k \in \mathcal{K}} \sum_{r \in \Omega_k} a_{i r}^k \iota^k=1, \quad \forall i \in \mathcal{C}\quad :\boldsymbol{\lambda}  \label{constr:mp_route1-lrp}\\
    & \sum_{k \in \mathcal{K}} \sum_{r \in \Omega_k}\sum_{j \in\mathcal{M}_i(t)} b_{j r}^k \iota^k\leq s_i,\quad :\boldsymbol{\mu}  \label{constr:mp_route-fcs-lrp} \\
    & \qquad \qquad \forall i \in\mathcal{F}, \forall t \in\mathcal{T}_i\nonumber\\
    & \sum_{k \in \mathcal{K}} \sum_{r \in \Omega_k} c_r^k \cdot \iota^k\leq \beta,\label{constr:mp-better-obj}\quad : \boldsymbol{\alpha} \\
    & \iota^k \in\{0,1\}, \forall r \in \Omega_{k}, \forall k \in \mathcal{K} \label{constr:mp_domain-lrp}\\
    & (\ref{upgrade_cons})-(\ref{y_cons}).
\end{align}
Similar to \eqref{CSP_obj}, the objective here also includes the annual equivalent investment cost \eqref{obj_mp0-invest} and the service cost \eqref{obj_mp0-energy}, with the latter expressed as cost associated with routes. The station charger capacity is incorporated as constraint \eqref{constr:mp_route-fcs-lrp}. 
Constraint \eqref{constr:mp-better-obj} corresponds to the requirement \eqref{better_obj} in \textbf{(P0)}. 

As the number of feasible routes is typically extremely large, recording a complete collection $\Omega_r^k$ may not be computationally feasible\footnote{We use the same notation $\Omega_r^k$ for both \textbf{(MP1)} and \textbf{(MP0)} for simplicity, but the two collections can be different.}. Hence, we will solve in an iterative fashion to extend the route set through column generation. Specifically, both \textbf{(MP1)} and \textbf{(MP0)} are replaced by their respective LP relaxations with a set of restricted routes $\Omega_k' \subset \Omega_k$ as well as $\iota^k \geq 0, r \in \Omega_k'$. We denote them as the restricted master problems \textbf{(RMP1)} and \textbf{(RMP0)}. A path generator for each type of vehicle $k$ is used to find additional columns (i.e. routes) by solving the following pricing problem, 
\begin{align}
    \min ~~ \varphi_n \qquad \mbox{s. to:} \quad (\ref{single_visit_cust})-(\ref{var_dom_fo}),
\end{align}
where $\varphi_n$ represents the reduced cost for \textbf{(RMP$n$)}. For \textbf{(RMP1)} and \textbf{(RMP0)}, the expressions for the reduced costs are respectively defined as
\begin{align}
    &\varphi_{1} = c_r^k-\sum_{i\in \mathcal{C}}a_{i r}^k\cdot\gamma_i,\\
    &\varphi_{0} = co_r^k-\sum_{i\in \mathcal{C}}a_{i r}^k\cdot \lambda_i+\sum_{i\in \mathcal{F}}\sum_{t\in \mathcal{T}_i}\mu_{i t}\cdot \sum_{j \in\mathcal{M}_i(t)} b_{j r}^k +\alpha\cdot c_r^k.
\end{align}
Note that $\gamma_i, \lambda_i, \mu_{it}, \alpha$ are dual variables associated with constraints \eqref{constr:mp_route1},\eqref{constr:mp_route1-lrp},\eqref{constr:mp_route-fcs-lrp},\eqref{constr:mp-better-obj}, respectively. This problem is also known as the Shortest Path Problem with Resource Constraints (SPPRC) \cite{irnich2005shortest, desaulniers2016exact, yu2019branch} and only paths with negative $\varphi$ are added into $\Omega_k'$.

The LP relaxations can be easily solved, whereas the SPPRC now accounts for the major computation burdens. A customized labeling algorithm is designed to solve the SPPRC.

\subsubsection{Labeling algorithm design} \label{sec:set_partitioning}
The main idea of a labeling algorithm is to find feasible paths in a dynamic programming fashion where the feasibility and superiority of each partial path are examined and compared along the process. To record the information required for this process, a label $L$ is defined for a partial path with an array of quantitative attributes. The specific label structure used in our work is defied as $L = \{node(L), w(L), \mathcal{R}(L), Path(L)\}$, where 

\begin{itemize}
    \item $node(L)$ is the last visited node of label $L$,
    \item $w(L)$ is the current cumulative cost of label $L$,
    \item $\mathcal{R}(L)$ is a set of the current cumulative resources consumed for each resource. Specifically, we have ${\mathcal{R}(L)=\left\{R_n(\cdot),R_d(\cdot), R_l(\cdot),R_e(\cdot),R_t(\cdot),R_{tw}(\cdot)\right\}}$  representing the number of nodes visited, total distance traveled, loaded capacity, the energy level and the departure time of the vehicle after serving node(L) and whether the time window constraint can be satisfied at current node,
    \item $Path(L)$ is the set of visited nodes in the traversed path up to $node(L)$.
\end{itemize}

Originating from $D_0$, a label is constructed iteratively by extending toward the possible descendant nodes in graph $\mathcal{G}^{PTE}$. A feasible extension needs to follow a set of resource criteria, namely resource extension functions (REFs). To extend a label $L$ from node $i$ to node $j$, the following REFs must be respected:
\begin{itemize}
    \item $0\leq R_n(i)+1\leq |\mathcal{N}_{0,n+1}|$,
    \item $0\leq R_d(i)+d_{ij}\leq D^\text{max}$,
    \item $0\leq R_l(i)+p_j\leq R_{l^\text{max}}$, 
    \item $0\leq R_e(i)+e_{ij}\leq R_{e^\text{max}}$, 
    \item $0\leq R_t(j) = R_t(i)+t_{ij}+t_{j}^{s} \leq T$, 
    \item $R_{tw}(j) =    
    \begin{cases}
    1 \quad \text{if } R_t(i)+t_{ij}\geq t_{i}^\text{min}, \text{and} \\ \qquad \quad R_t(i)+t_{ij}+t_{j}^{s}\leq t_{j}^\text{max}, \\
    0 \quad \text{otherwise}.
    \end{cases}$
\end{itemize}
\subsubsection{Acceleration Techniques}
There are two major techniques we adopt and apply to accelerate the inner-loop process. Namely, they are stabilized column generation and dominance check. These two techniques are particularly useful for reducing the number of necessary iterations and accelerating individual iterations, respectively.

\textbf{Stabilized Column Generation: }We perturb the problem at constraint \eqref{constr:mp_route1} with bounded surplus and slack variables to overcome primal degeneracy such that a trust region is formed around the associated dual variables. Notice that these duals play an important role in guiding the exploration of feasible routes. The trust region is helpful in obtaining an accurate dual estimation and hence reduces the number of needed iterations. This is the well-known BOXSTEP method and the perturbed version of the problem has been shown to be equivalent to the original. We refer the readers to \cite{du1999stabilized,lubbecke2005selected} for further details.

\textbf{Dominance Check in Labeling Algorithm: }It should be noticed that it is computationally inefficient to exhaust extensions for every single label. Hence, the dominance check is implemented to recognize and eliminate partial paths that need no further examination, since they are dominated by other partial paths.
A label $L_1$ is said to dominate label $L_2$ if the followings dominance rules are all satisfied:
\begin{itemize}
    \item $node(L_1)$ = $node(L_2)$,
    \item $w(L_1) \leq w(L_2)$,
    \item $R(L_1) \preceq R(L_2)$,
    \item $Path(L_2)\subseteq Path(L_1)$,
\end{itemize}
where $\preceq$ is used for element-wise comparison. By removing all labels like $L_2$, the number of potential extensions can be largely reduced. The pseudo code for the overall labeling algorithm is presented in Algorithm \ref{algo:labeling}.

 \begin{algorithm}[ht!]
 \caption{Labeling Algorithm}
 \begin{algorithmic}[1]
 \renewcommand{\algorithmicrequire}{\textbf{Input:}}
 \renewcommand{\algorithmicensure}{\textbf{Output:}}
 \REQUIRE network and vehicle specification; REFs
 \ENSURE optimal path for freight of type $k$
 \\ \textbf{\textit{- Initialization}} :
  \STATE \textit{[max\_res], [min\_res]: } \textit{UB} and \textit{LB} arrays of resource usage for freight of type $k$; 
 \\ \textit{L\_current: } \{0, Source, \textit{[min\_res]}, [``Source"]\};
 \\ \textit{\{unprocessed\_labels (unp\_labels)\}: } empty queue.
 \\ \textbf{\textit{- Loop Process :}}
  \WHILE {\textit{L\_current} != None}
  \FOR{edge in $\mathcal{G}$.$node(\textit{L\_current})$.adjacent\_edges}
  \STATE PROPAGATE\_LABEL(edge, \textit{L\_current}): \textit{L\_edge}
  \IF{\textit{L\_edge} feasible}
  \STATE Add to \textit{\{unp\_labels\}}
  \ENDIF
  \ENDFOR
  \STATE EXTEND\_LABEL(\textit{L\_current}, \textit{\{unp\_labels\}}): \textit{L\_next}
      \STATE Update \textit{L\_current} = \textit{L\_next}
  \STATE CHECK\_DOMINANCE(\textit{L\_current}, \textit{\{unp\_labels\}})
  \IF{\textit{Label} being dominated}
  \STATE Delete from \textit{\{unpr\_labels\}}
  \ENDIF
  \STATE SAVE\_CURRENT\_BEST\_LABEL(\textit{L\_current}): \textit{L\_final}
  \ENDWHILE
 \STATE \textit{L\_best} = \textit{L\_final}
 \RETURN \textit{L\_best}
 \end{algorithmic}
 \label{algo:labeling}
 \end{algorithm}
\section{Case Studies} \label{sec:case_study}
We have proposed an optimization modeling framework and iterative algorithm to solve this problem. To effectively demonstrate the model, we deliberately design a small but intuitive network. We will highlight some of the binding features, like the customer time windows and the charging rates. 

\begin{figure}[!tp]
  \begin{center}
  \includegraphics[scale=0.6]{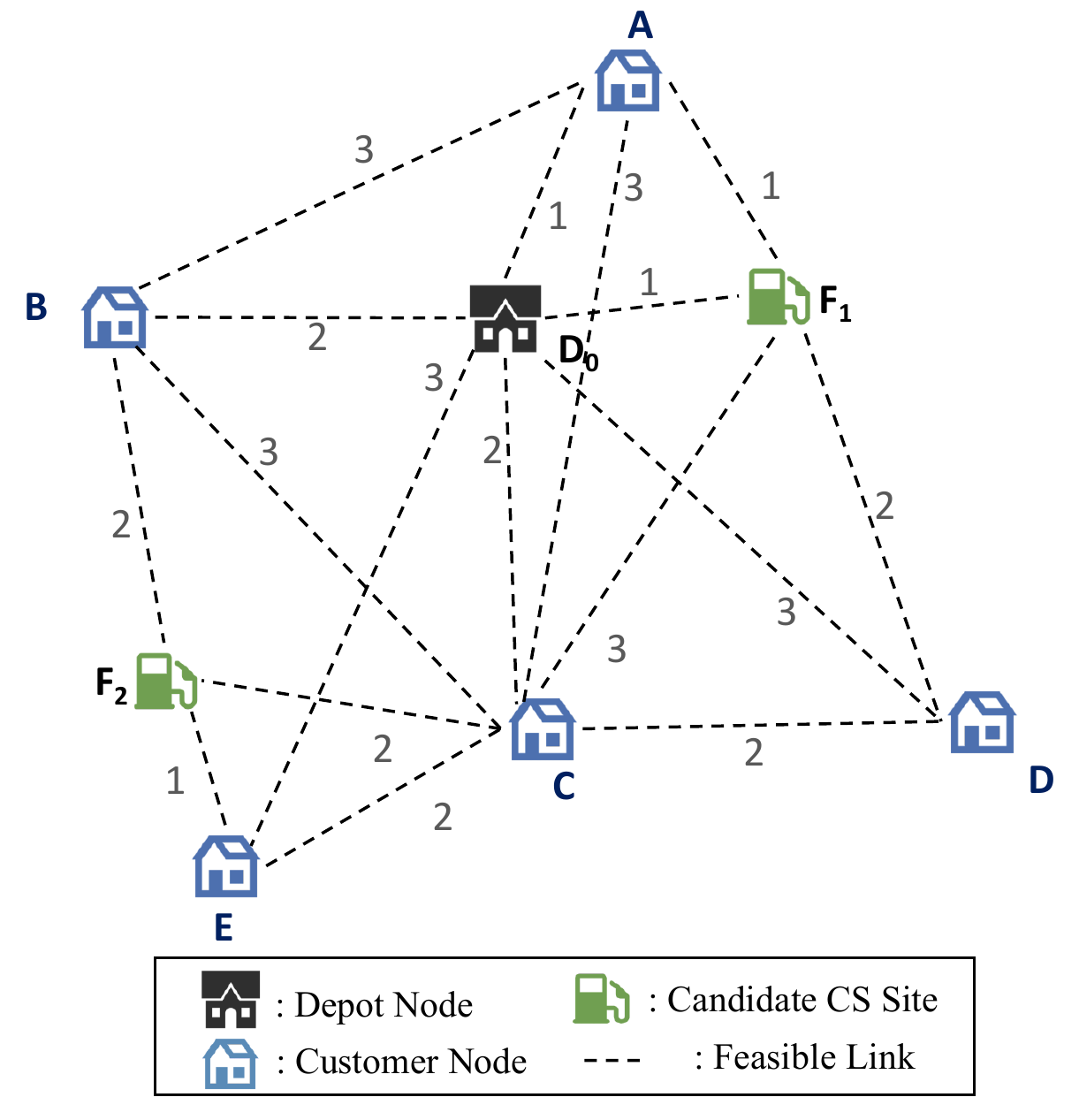}
  \caption{Small network (original)}\label{fig-casestudy_net}
  \end{center}
\end{figure}

The small network is presented in Fig. \ref{fig-casestudy_net}. It consists of 1 depot node $D_0$, 5 customer nodes $\mathcal{C} = \{A, B, C, D, E\}$ and 2 candidate charging station nodes $\mathcal{F} = \{F_1, F_2\}$. The dashed lines are the feasible links with adjacent numbers indicating the lengths. We summarize all relevant parameters for CSP in Table \ref{table:small-case-FCS} and for the FO in Table \ref{table:small-case-vehicle}. Next, we will numerically demonstrate the necessity of the two-entity modeling. Then, we will study the cost breakdowns and the varying dynamics when stricter time windows are applied and charging rates are varied. 

\begin{table*}[!t]
\renewcommand{\arraystretch}{1.3}
\caption{Parameters of candidate charging stations}
\label{table:small-case-FCS}
\centering
\begin{tabular}{|c|c|c|c|c|c|}
\hline
\makecell{Station\\ID} & \makecell{Available\\capacity (kW)}& \makecell{Charge\\rate (kW)} & \makecell{Cost of\\charger (\$)} & \makecell{Substation\\ upgrade cost (\$/kW)} & \makecell{Electricity\\cost (\$/kWh)}\\
\hline
$F_1$ & 15 & 5	& 10000 & 788 & 0.1\\
\hline
$F_2$ & 30 & 10 & 10500 & 788 & 0.1\\
\hline
\end{tabular}
\end{table*}

\begin{table*}[!t]
\renewcommand{\arraystretch}{1.3}
\caption{Vehicle parameters}
\label{table:small-case-vehicle}
\centering
\begin{tabular}{|c|c|c|c|c|c|}
\hline
\makecell{Vehicle\\ID} & \makecell{Freight\\capacity (kg)} & \makecell{Battery\\capacity (kWh)}& \makecell{Energy\\consumption (kW/unit length)} & \makecell{Vehicle\\cost (\$)} & \makecell{Travel\\cost (\$/unit length)} \\
\hline
1 & 150 & 50	& 10 & 10000 & 0.5\\
\hline
2 & 200 & 60 & 10 & 18000 & 0.5 \\
\hline
\end{tabular}
\end{table*}

\subsection{Base Case: Necessity of Considering Different Entities} \label{sec:base_case}
\begin{figure}[!tp]
  \begin{center}
  \includegraphics[width=\columnwidth]{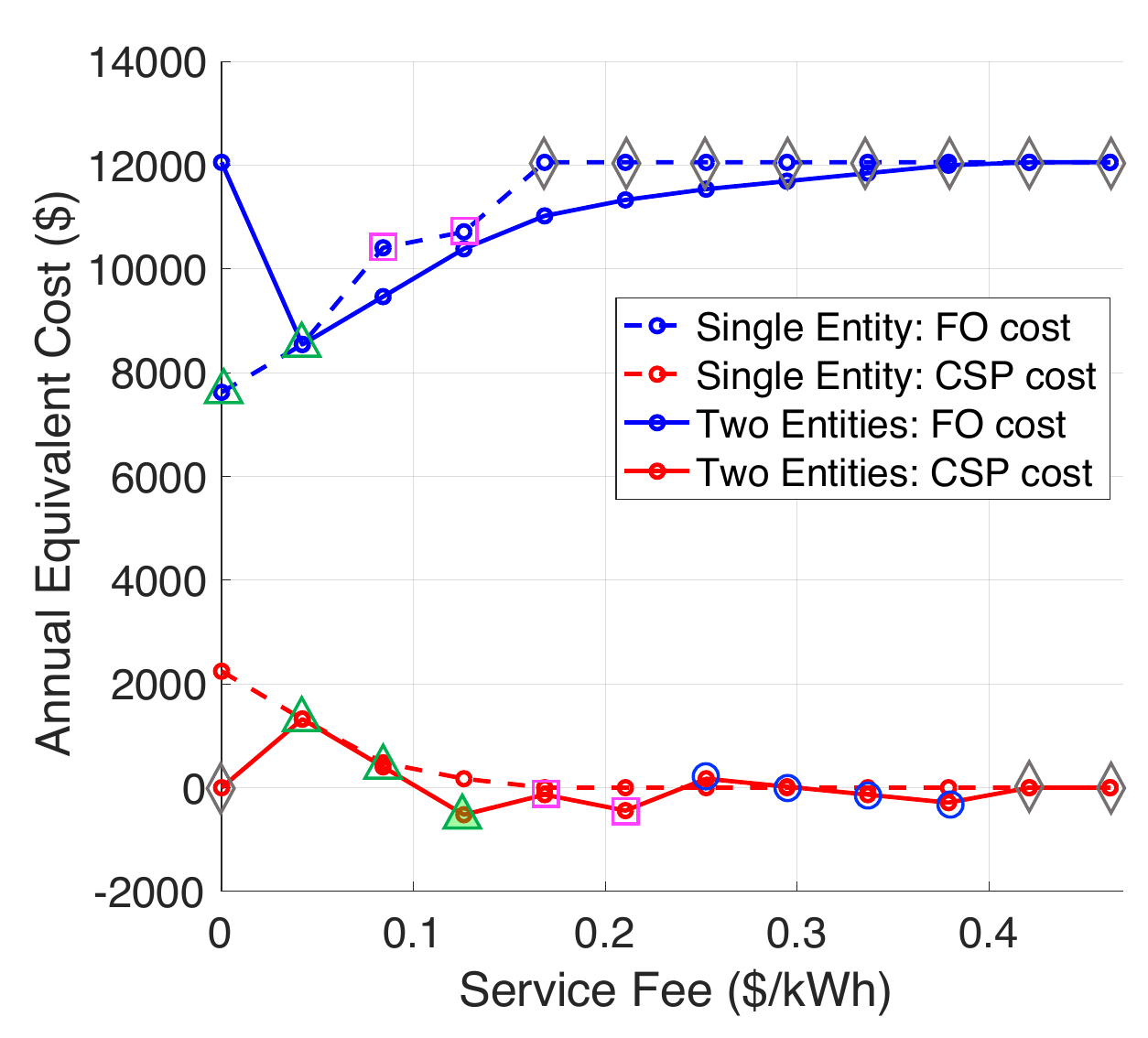}
  \caption{Cost of different players with respect to different service fee, without time windows. The colored marker shapes correspond to different strategies in Table VI. The markers on the upper dashed blue line correspond to the single entity case, whereas the markers on the lower solid red line correspond to the two-entity case.}\label{fig:Rate10_without_TW-lowSingle}
  \end{center}
\end{figure}

We first focus on the different strategies when a social planner (single entity) or non co-operation (two entities) is considered. Simulations are performed for both scenarios with the service fee varying from 0 to 0.5($\$$/kWh) \footnote{We leave customer time windows sufficiently wide in this case.}. The results obtained are given in Fig.\ref{fig:Rate10_without_TW-lowSingle}.

Consider the single entity scenario. To enable a fair comparison against the two-entity case, we plot the combined costs of electricity and service fee ($\pi_i+c_i^e$) for the single entity. Together they will jointly affect the vehicle routing, charging, as well as the infrastructure decisions. When the service fee is set to 0, the FO charges at the cost of electricity ($\pi_i$) purchased from the utility, i.e. with zero profit margin. Given the optimized planning results, we then split and plot the corresponding costs to the FO and CSP. In greater detail, the FO cost consists of the fleet investment, travel expenditure as well as the combined cost of electricity and service fee; on the other hand, the CSP cost is the infrastructure investment less the profits from providing service. The cost splits are presented by the dashed lines in Fig. \ref{fig:Rate10_without_TW-lowSingle} and the colored shape labels on the dashed line indicate the different optimal strategies corresponding (see Table VI) to each simulated price value. Note that when service fee is set at \$0/kWh, this is the case commonly known as to minimize the total cost of ownership (TCO) in literature. However, as shown in the left most in Fig.\ref{fig:Rate10_without_TW-lowSingle}, it actually induces the largest cost to the CSP.


For the two-entities scenario, we solve the problem with the proposed model. We superimpose the costs for the FO and CSP in Fig. \ref{fig:Rate10_without_TW-lowSingle} using solid lines. It is visually clear that the decisions under the two-entity scenario achieves lower net costs than the single entity scenario most of the time, i.e. the FO achieves lower costs and the CSP generates more profit. The only exceptions happen when the service fees are extremely low (below \$0.05/kWh). Within this range, from a pure economic point of view, the CSP has little to no interest to invest and enter the market \footnote{Even if the CSP enters (the second green triangle), it is not at all cost-attractive. It experiences positive cost and no profit.}. The CSP can be incentivized to participate by increasing the service fee. However, as the CSP becomes more ``greedy,'' the FO will reject the option to charge at the facility. This then leads to a non-cooperative situation, as indicated by the grey rhombuses on the solid line in Fig. \ref{fig:Rate10_without_TW-lowSingle} (after service fee is more than \$0.421/kWh).

\begin{figure}[!tp]
  \begin{center}
  \includegraphics[scale=1]{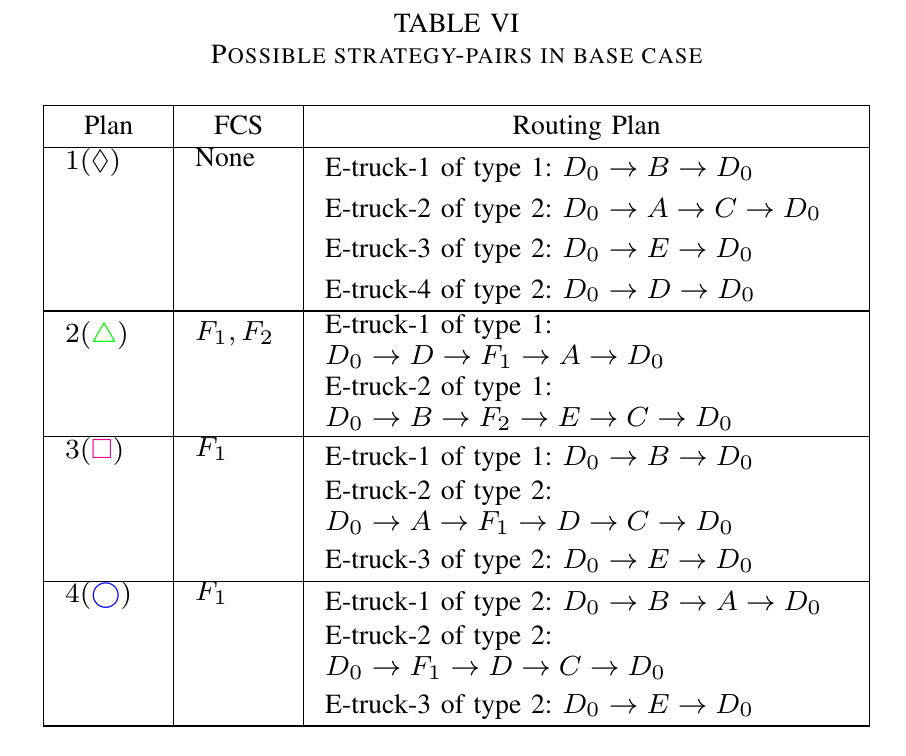}
  \label{table:small-case-strategyPairs}
  \end{center}
\end{figure}
\begin{figure}[!tp]
  \begin{center}
  \includegraphics[scale=0.6]{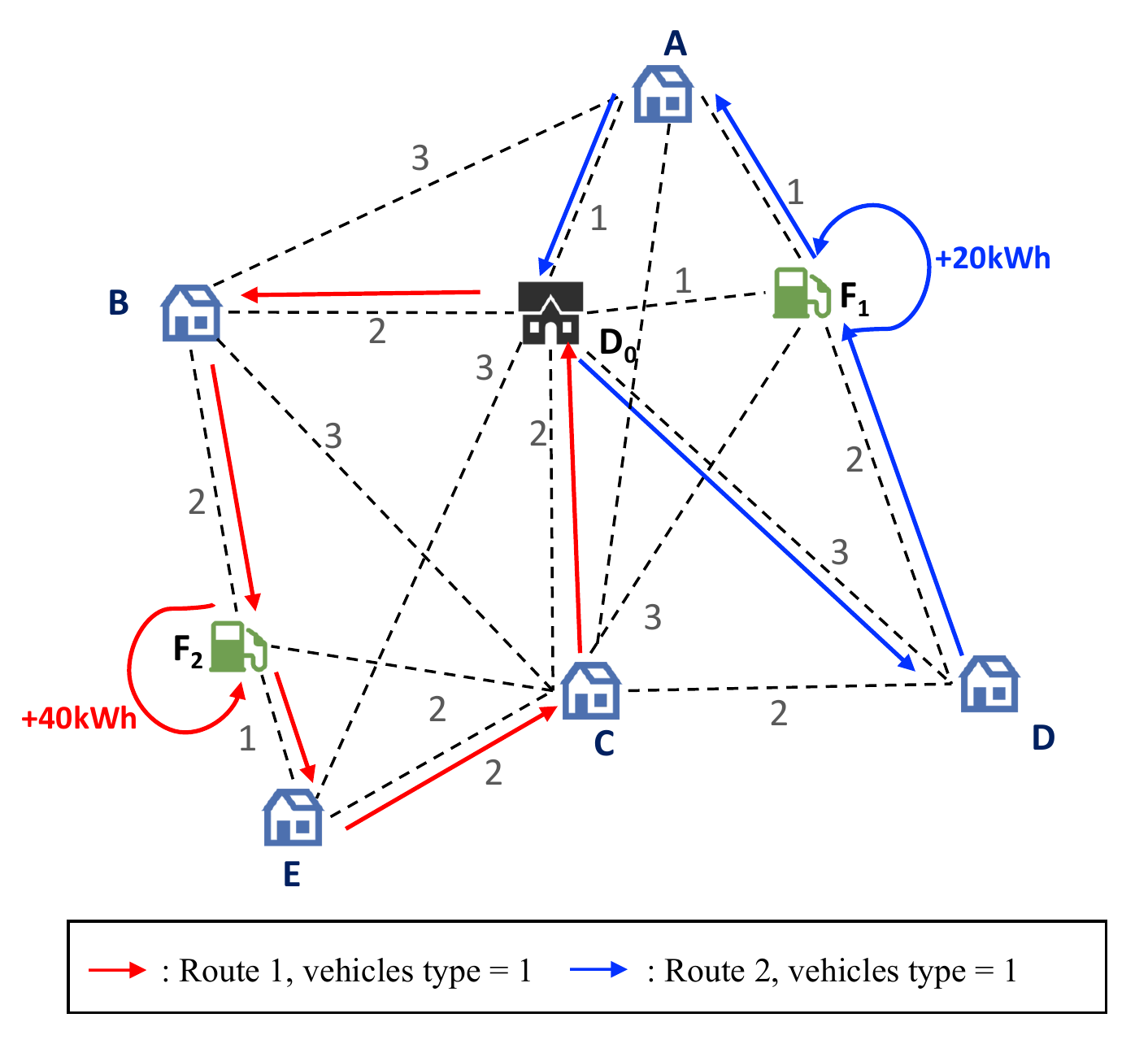}
  \caption{Optimal routes of FO. }\label{fig:Strategy_TwoEntity_noTW}
  \end{center}
\end{figure}

Since the CSP is the leader in the game, it will naturally take advantage of being the first-mover. The final decision for the CSP is to build charging stations at both locations with service fee set at 0.125 \$/kWh in this particular example (solid green triangle in Fig. \ref{fig:Rate10_without_TW-lowSingle}). The corresponding fleet routes are shown in Fig. \ref{fig:Strategy_TwoEntity_noTW} (Plan 2). We denote the first/second route chosen by E-truck-1/2 as route 0/1 for later reference. 

Indeed, it is the most efficient and straightforward to assume the existence of a powerful social planner. However, this common assumption in the literature is not always the case in real life. We further emphasize an important nuance here: only through rigorous and comprehensive cost analysis can we incentivize both entities to come up with an agreement that benefits both sides, even though the eventual outcome does not always lead to the social optimum as contests prevail \footnote{Note that this is an one-to-one service scenario. More involved settings and further analysis can be extended to one-to-many or many-to-one scenarios. This is out of scope and we will leave it for future work.}.


\subsection{The effects of customer time windows and charging rates}
\begin{figure}[!tp]
  \begin{center}
  \includegraphics[width=\columnwidth]{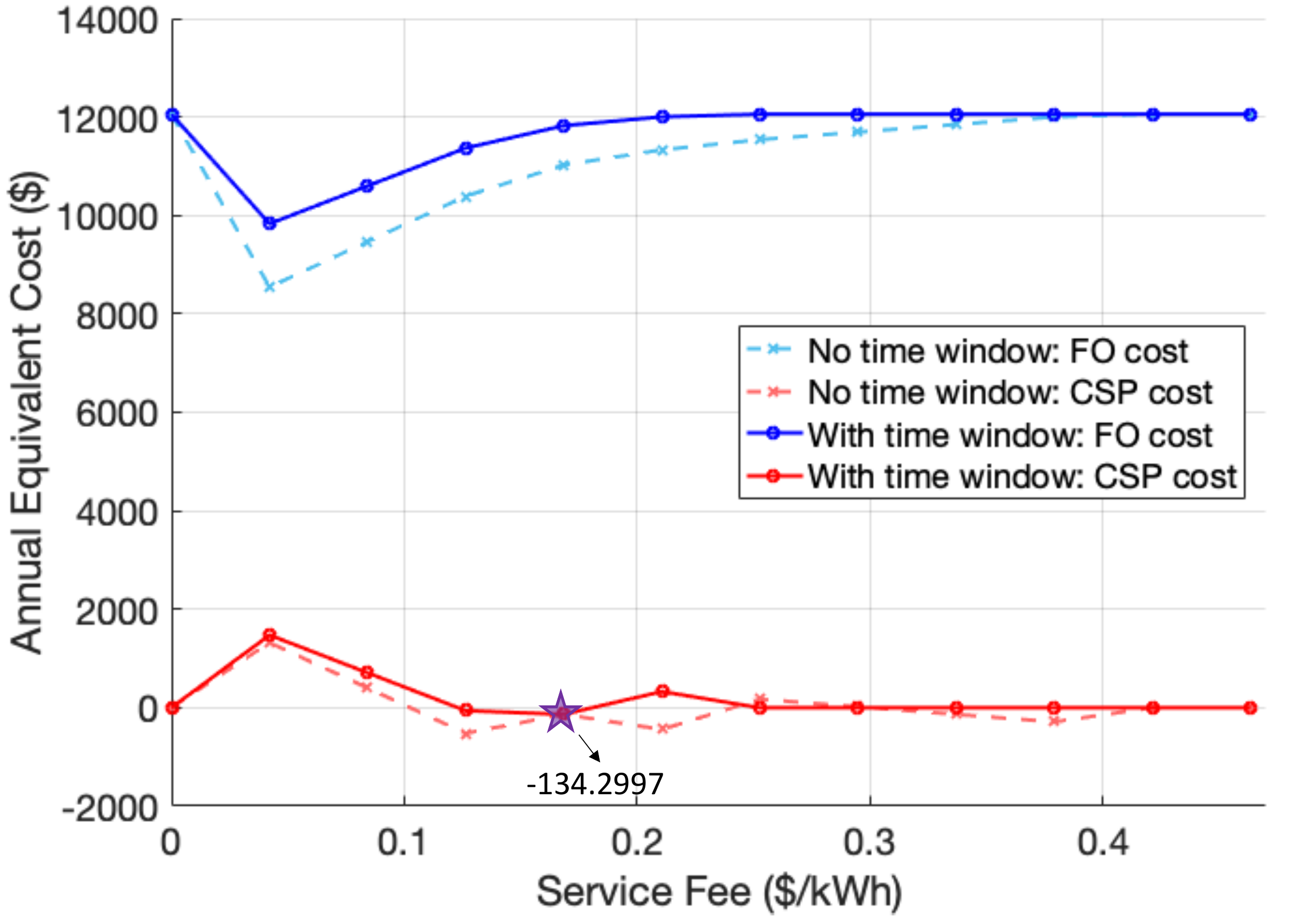}
  \caption{With time windows, cost of different players with respect to different service fee. The two upper curves correspond to the FO costs; the lower two are the CSPs'.}\label{fig-Rate10_with_TW}
  \end{center}
\end{figure}

We are going to show how two of the binding features, time windows and charging rates, will influence optimal decisions. We set tighter delivery time windows for customer $B, C$ and $D$ as $[1,4],[6,9],[2,3]$ respectively. First, we keep the charge rate as 10 kW as in the base case and the results are plotted in Fig. \ref{fig-Rate10_with_TW}. In this case, the FO will reject charging service from the CSP when the price exceed $0.25$ \$/kWh. This critical threshold was $0.421$ \$/kWh in the base case (Section \ref{sec:base_case}), where no time windows are set. The only profitable service fee for the CSP occurs $0.1684$ \$/kWh, whereas all other services fees yield zero or negative profit. Note, in the base case, there were six profitable service fees. However, the profit now is $74.23\%$ lower than that in the base case. 

\begin{figure}[!tp]
  \begin{center}
  \includegraphics[scale=1]{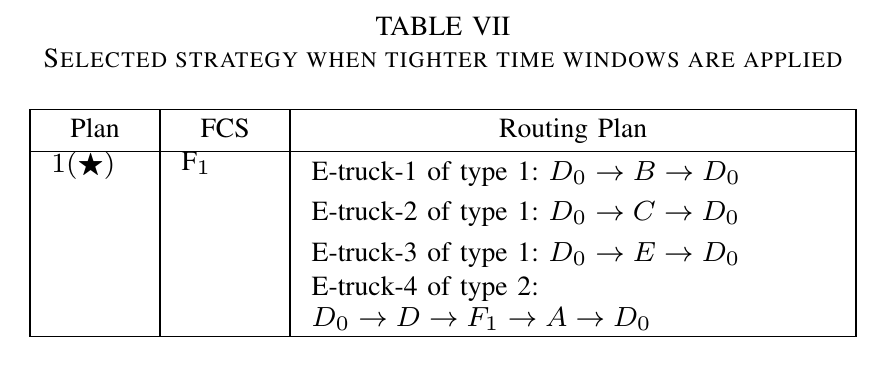}
  \label{table:small-case-withTW_strategy}
  \end{center}
\end{figure}
In this example, some tightened time windows significantly reduced the number of feasible routes for the FO. It becomes a tougher decision to balance both time spent for charging and time to meet customers. The previously selected route 1 ($D_0 \rightarrow B\rightarrow F_2\rightarrow E \rightarrow C \rightarrow D_0$) becomes infeasible as the E-trucks need 4 units of time to charge (to complete the trip) but will miss customer $C$. As a result, three separate E-trucks are purchased instead to serve customer $B, C$ and $E$ respectively. The detailed plan is given in Table VII. This increases the final cost of the FO by $13.8\%$.

\begin{figure}[!tp]
  \begin{center}
  \includegraphics[width=\columnwidth]{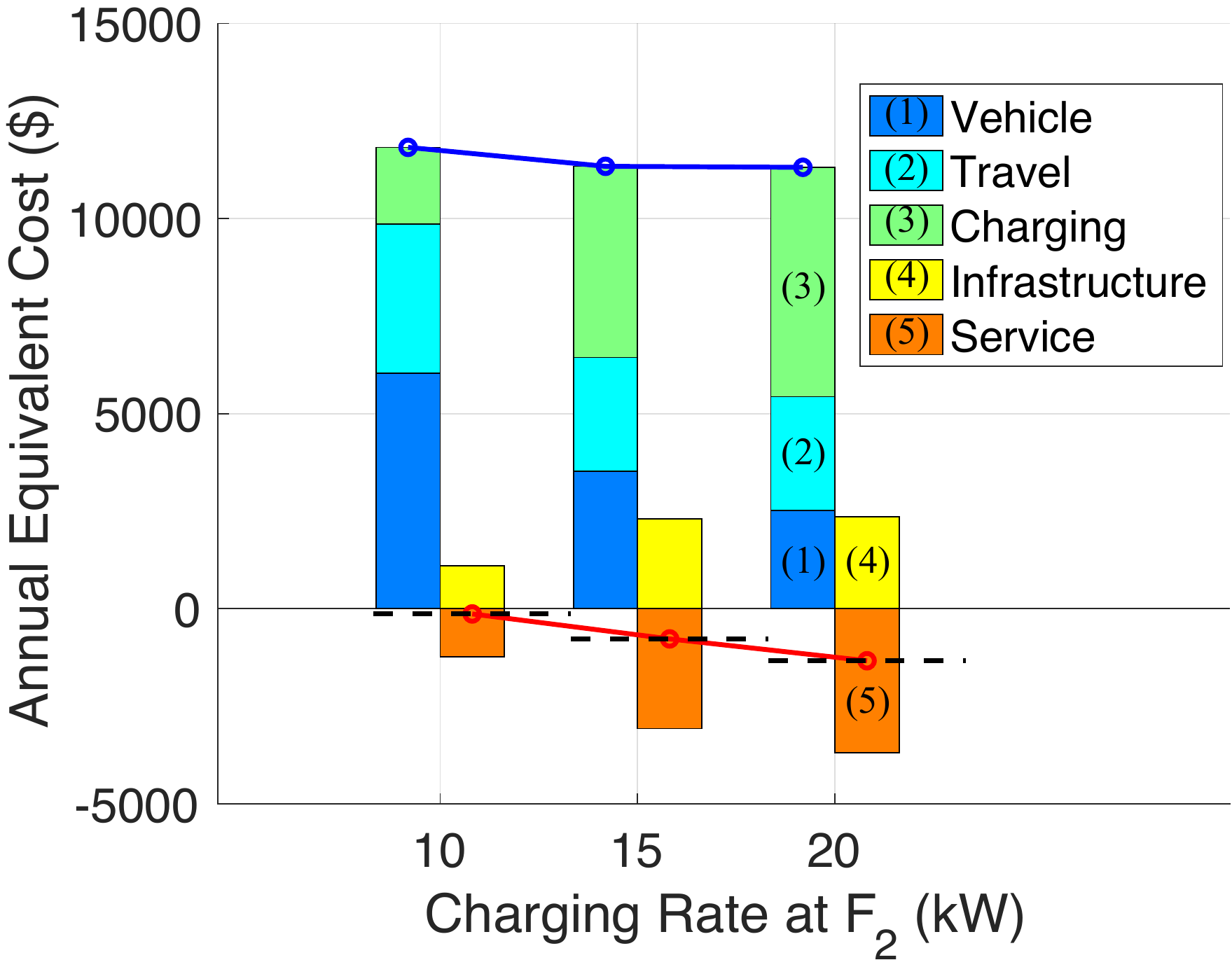}
  \caption{Costs of different players with different charging rate at $F_2$. Upper solid line is the net annual equivalent costs of the FO and the lower dashed line is the net profits of the CSP.} \label{fig:cost_breakdown_wrtRate}
  \end{center}
\end{figure}

One alternative solution that may realize benefits to both sides is to increase the chargers' charging rate. The cost breakdowns over the variations of rate are shown in Fig. \ref{fig:cost_breakdown_wrtRate}. They are plotted with respect to the two entities -- the FO on the left and the CSP on the right. Black dashed lines indicate the CSP's overall net cost/profit (for positive/negative values resp.). The corresponding strategies are also listed in Table VIII. We see that when the rate increases from 10 kW to 15 kW, the previously selected route 1 can again be assigned to a larger E-truck, whereas when the rate is doubled, the smaller E-truck may be used. Although more charging energy is needed, the trucks still meet all the delivery time windows. As a result, reductions in E-truck fleet investment and travel cost compensate the extra charging expenditure and lead to overall FO cost savings. From the CSP's perspective, although more infrastructure investment is required\footnote{The cost for a single port is assumed to increase with the charging rate. Indeed the infrastructure investments at 15 and 20 kW are subtly different.} up-front, it also realizes more profit gains. The increased revenues from energy service cover the increased capital costs. 



\begin{figure}[!tp]
  \begin{center}
  \includegraphics[scale=1]{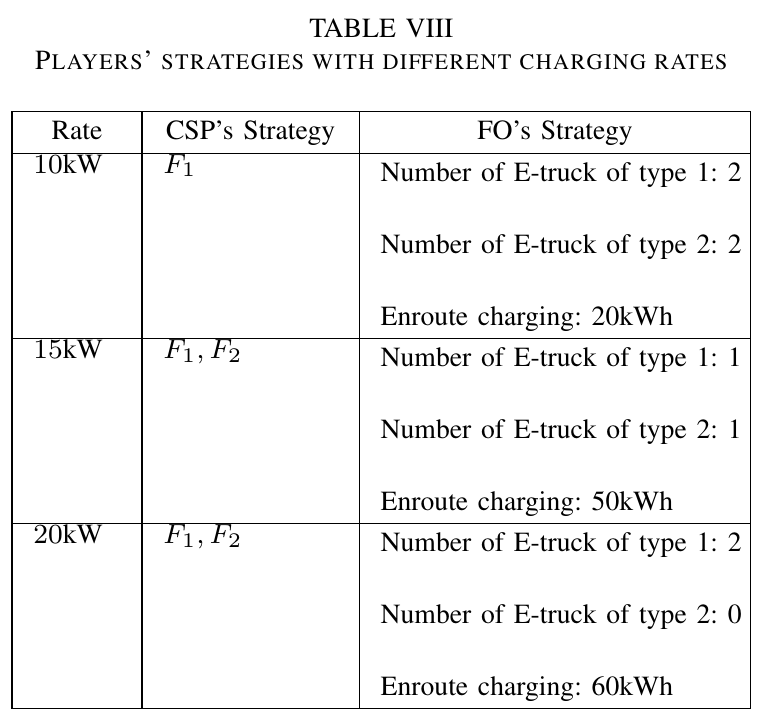}
  \label{table:strategy-with-diffRate}
  \end{center}
\end{figure}

We have illustrated two of the binding features along with the dynamic interactions between the two entities. This is indeed a complex but more realistic situation, as opposed to existing work that considers a single entity only. We have demonstrated that without assuming the existence of a powerful social planner, careful analysis is needed to contemplate and incentivize both parties to collaborate. With our framework, the FOs and the CSPs can now look for their operation ``sweet spots'' and be more committed to transportation electrification.

\section{Conclusion} \label{sec:conclusion}
This work proposes a new modeling framework to capture the non-cooperative interactions between a charging service provider (CSP) and fleet operator (FO) in the joint charging network design and mobility planning problem. This reflects reality, in which the two entities are separate self-optimizing organizations. To consider the charging station capacity, this work also proposes, for the first time, a partial time expanded network. It enables jointly optimizing the size of the charging infrastructure in the classic location routing problem regime. The non-cooperative interaction is formulated in a Stackelberg game framework as a Bi-MILP problem. A double-loop solution framework with customized labeling algorithm is designed accordingly.

We find that non-cooperation between these two entities can lead to lower total costs than a single entity scenario. Detailed numerical studies demonstrate the effects of binding features and provide insights to interested players. The nature of the solutions with a single social planner and two non-cooperative entities can be dramatically different. 
Our framework examines these distinct aspects and hence provide more reliable results to the case in real life.
In fact, multiple entities could co-exist in many cases. A future research direction may also consider the one-to-multiple entity scenario.

\section{Acknowledgement}
This work was supported in part by Total R\&D. The authors would like to thank Dr. Anna Robert for providing constructive feedbacks.


%





\ifCLASSOPTIONcaptionsoff
  \newpage
\fi





\bibliographystyle{IEEEtran}
\bibliography{IEEEabrv,Bibliography}

\vfill


\end{document}